\documentclass[journal]{IEEEtran}

\usepackage{tikz}
\graphicspath{{images/}}
\usepackage{amsmath}
\usepackage{latexsym, amssymb}
\usepackage{graphicx}
\usepackage{amsmath, amsbsy}
\usepackage{amsopn, amstext}
\usepackage{ifpdf,hyperref}
\usepackage{cancel, color}
\usepackage{epstopdf}

\def\lplus{\,\rotatebox[]{-90}{$\pm$}\,}

\def\lminus{\vdash}

\def\lvtimes{\vec{\ltimes}}

\def\lvminus{\vec{\vdash}}

\def\lvplus{\vec{\lplus}}

\def\J{{\bf 1}}

\def\A{\left< A\right>}
\def\B{\left< B\right>}
\def\C{\left< C\right>}
\DeclareMathOperator{\bed}{bed}

\DeclareMathOperator{\Span}{Span}
\DeclareMathOperator{\Col}{Col}
\DeclareMathOperator{\Row}{Row}

\def\cal{\mathcal}

\def\ra{\rightarrow}

\def\lra{\leftrightarrow}

\def\a{\alpha}
\def\b{\beta}
\def\d{\delta}

\def\0{{\bf 0}}

\newcommand{\R}{{\mathbb R}}
\newcommand{\Q}{{\mathbb Q}}
\newcommand{\N}{{\mathbb N}}

\def\dsum{\mathop{\sum}\limits}

\newtheorem{thm}{Theorem}[section]
\newtheorem{dfn}[thm]{Definition}
\newtheorem{prp}[thm]{Proposition}
\newtheorem{exa}[thm]{Example}
\newtheorem{lem}[thm]{Lemma}
\newtheorem{cor}[thm]{Corollary}
\newtheorem{rem}[thm]{Remark}

\begin{document}

\title{Cross-Dimensional Linear Systems}

\author{Daizhan Cheng, Zequn Liu, Hongsheng Qi
	\thanks{This work is supported partly by the National Natural Science Foundation of China (NSFC) under Grants 61733018, 61333001, and 61773371.}
	\thanks{Daizhan Cheng is with the Key Laboratory of Systems and Control, Academy of Mathematics and Systems Sciences, Chinese Academy of Sciences,
		Beijing 100190, P. R. China (e-mail: dcheng@iss.ac.cn).}
    \thanks{Zequn Liu and Hongsheng Qi are with the Key Laboratory of Systems and Control, Academy of Mathematics and Systems Sciences, Chinese Academy of Sciences, Beijing 100190, and University of Chinese Academy of Sciences, Beijing 100049, P. R. China (e-mail: liuzequn@amss.ac.cn, qihongsh@amss.ac.cn).}
    \thanks{Corresponding author: Daizhan Cheng. Tel.: +86 10 82541232.}
}

%\markboth{IEEE Transactions On Automatic Control, Vol. XX, No. Y, Month 201Z}
%{Cheng \MakeLowercase{\textit{et al.}}: On Dimension-Free Linear System}

\maketitle

\begin{abstract}
Semi-tensor product(STP) or matrix (M-) product of matrices turns the set of matrices with arbitrary dimensions into a monoid $({\cal M},\ltimes)$. A matrix (M-) addition is defined over subsets of a partition of ${\cal M}$, and a matrix (M-) equivalence is proposed. Eventually, some quotient spaces are obtained as vector spaces of matrices. Furthermore, a set of formal polynomials is constructed, which makes the quotient space of  $({\cal M},\ltimes)$, denoted by $(\Sigma, \ltimes)$, a vector space and a  monoid.

Similarly, a vector addition (V-addition) and a vector equivalence (V-equivalence) are defined on ${\cal V}$, the set of vectors of arbitrary dimensions. Then the quotient space of vectors, $\Omega$, is also obtained as a vector space.

The action of  monoid  $({\cal M},\ltimes)$ on ${\cal V}$ (or $(\Sigma, \ltimes)$ on $\Omega$) is defined as a vector (V-) product, which becomes a pseudo-dynamic system, called the cross-dimensional linear system (CDLS). Both the discrete time and the continuous time CDLSs  have been investigated. For certain time-invariant case, the solutions (trajectories) are presented. Furthermore, the corresponding cross-dimensional linear control systems are also proposed and the controllability and observability are discussed.

Both M-product and V-product are generalizations of the conventional matrix product, that is, when the dimension matching condition required by the conventional matrix product is satisfied they coincide with the conventional matrix product. Both M-addition and V-addition are generalizations of conventional matrix addition. Hence, the dynamics discussed in this paper is a generalization of conventional linear system theory.

\end{abstract}

\begin{IEEEkeywords}
Monoid, dynamic-system, cross-dimensional linear system, M-product/M-addition of matrices, V-product/V-addition of matrices/vectors.
\end{IEEEkeywords}

\IEEEpeerreviewmaketitle

\section{Preliminaries}%s-1

Dimension-varying system exists everywhere. For instance, in the internet or some other service-based networks, some users may join in or withdraw out from time to time. In a genetic regulatory network, cells may die or birth at any time. To the authors' best knowledge, so far there is no proper tool to handle or even model dimension-varying system properly.

From engineering point of view, a single electric power generator may be described by 2-dimensional model, or 3-dimensional, or even 5-, 6-, or 7-dimensional models \cite{mac97}. From theoretical physics, the superstring theory proposes the space-time dimension could be 4 (Instain, Relativity), 5(Kalabi-Klein theory), 10 (Type 1 string), 11 (M-theory), or even 26 (Bosonic), etc. \cite{kak99}

It seems that we may need a theory to connect spaces with different dimensions together. In other words, we may need to allow a dynamic system to go cross spaces of different dimensions.
The purpose of this paper is to explore a cross dimension linear system. Though the proposed model may help to solve the modeling and analysis of dimension-varying systems, the work in this paper is mainly of mathematical interest.  There is a long way ahead to go before practically solving the modeling, analysis and control of dimension-varying systems. The break through point in this work is: the distance between two points on spaces of different dimensions separately can be measured. Then the difference between two systems on different dimension spaces can also be measured.

For a CDLS, the dimension of state variables could be arbitrary and time-varying. Precisely speaking, choosing any initial state $x_0\in \R^s$, where $s\in \N$ is arbitrary, then the system is well defined and evolving over dimension-varying state space.

The theoretical starting point of this approach is the semi-tensor product (STP) of matrices, which is also called the matrix (M-) product in this paper. The STP is a generalization of conventional matrix product and it can be used for two arbitrary matrices \cite{che11,che12}. In the light of M-product, the matrices of arbitrary dimensions become a monoid. The CDLS is based on the action of monoid of matrices with M-product on dimension-varying vector space as the state space. The monoid-based dynamic system is briefly introduced as follows.

\begin{dfn}\label{d1.1} Let $G$ be a nonempty set, and $*$ a binary operator on $G$.
\begin{enumerate}
\item $(G,*)$ is called a semi-group, if (Associativity)
\begin{align}\label{1.1}
(g_1*g_2)*g_3 =g_1*(g_2*g_3),\quad g_1,g_2,g_3\in G.
\end{align}
\item A semi-group $(G,*)$ with identity $e\in G$ is called a monoid. The identity satisfies
\begin{align}\label{1.2}
e*g =g*e=g,\quad g\in G.
\end{align}
\end{enumerate}
\end{dfn}

\begin{rem}\label{r1.101} A monoid becomes a group if for each $g\in G$ there exists a $g^{-1}$ such that
\begin{align}\label{1.201}
g*g^{-1} =g^{-1}*g=e,\quad g\in G.
\end{align}
\end{rem}

\begin{dfn}[\cite{liu08,how95}]\label{d1.2} Let $S$ be a monoid and $X\neq \emptyset$ an objective state space. A mapping $\varphi:S\times X\ra X$ is called an S-system of $S$ on X, if it satisfies
\begin{itemize}
\item
\begin{align}\label{1.3}
\varphi(g_1,\varphi(g_2,x))=\varphi(g_1*g_2,x),\quad g_1,g_2\in G,\;x\in X.
\end{align}
\item
\begin{align}\label{1.4}
\varphi(e,x)=x,\quad x\in X.
\end{align}
\end{itemize}
\end{dfn}

Hereafter, an S-system is denoted by $(S,\varphi,X)$, where $S$ is a monoid, $X$ is the objective state space, and $\varphi:~S\times X\ra X$ is said to be the action of $S$ on $X$. The following definition is similar to the one in \cite{kop75}.

\begin{dfn}\label{d1.3} An S-system $(S,\varphi,X)$  is called a pseudo dynamic system, if $X$ is a topological space and for each $s\in S$, $\varphi|_s:~X\ra X$ is continuous. If $X$ is a Hausdorff space, the pseudo dynamic system  becomes a dynamic system.
\end{dfn}
For notational ease, denote the action as:
$$
\varphi(s,x):=sx.
$$
Then an  S-system can be expressed in a conventional form as
\begin{align}\label{1.5}
x(t+1)=A(t)x(t),\quad x(t)\in X,\; A(t)\in S.
\end{align}

Roughly speaking, a cross dimension linear system considered in this paper is of two categories:
\begin{itemize}
\item  An S-system which is constructed by  (i) a monoid consisting of all matrices with M-product as its operator; and (ii) a pseudo-vector space of arbitrary dimensions. The action of prier on latter is the V-product.
\item An dynamic system which is constructed by (i) a quotient monoid of equivalent classes of matrices, as a vector space, with M-product as its operator; and (ii) a vector space of equivalent classes of vectors. The action of prier on latter is the V-product.
\end{itemize}

The STP (or M-product) of matrices was proposed about two decades ago.  After almost two decades, the STP has been achieved various applications for analysis and control of logical systems \cite{lupr,las13}, finite games \cite{guo13,che15}, fault detection \cite{for15}, graph theory \cite{wan12}, etc. Recently, the mathematical annotation of STP has been dug in by a long paper \cite{chepr}, which provides a mathematical background for this investigation.

The rest of this paper is organized as follows: Section 2 provides a framework for the S-system of matrix monoid on vector space of arbitrary dimensions. Section 3 considers the discrete time CDLS. Section 4 proposes formal polynomial, which has matrices of varying dimensions as its coefficients. Based on formal polynomial, Section 5 investigates the continuous time CDLS. The controlled CDLS is discussed in Section 6. Section 7 is a brief conclusion.

Before ending this section we give some notations:

\begin{enumerate}
\item $\N$: Set of natural numbers (i.e., $\N=\{1,2,\cdots\}$);

\item $\Q$: Set of rational numbers, ($\Q_+$: Set of irreducible positive rational numbers.);

\item $\R$: Field of real numbers;
%
%
%%\item $\F$: certain field of characteristic number $0$ (Particularly, we can understand $\F=\R$ or $\F=\C$).
%
\item ${\cal M}_{m\times n}$: set of $m\times n$ dimensional real matrices.

\item $\Col(A)$ ($\Row(A)$): the set of columns (rows) of ~$A$; $\Col_i(A)$ ($\Row_i(A)$): the $i$-th column (row) of ~$A$.

\item One vector: ${\bf 1}_n=[\underbrace{1,\cdots,1}_n]^T$.
\item
$a|b$: $a$ divides $b$.

\item
$m\wedge n=gcd(m,n)$: The greatest common divisor of $m$ and $n$.

\item
$m\vee n=lcm(m,n)$: The least common multiple of $m$ and $n$.

\item
$\left<X,Y\right>$, $X,Y\in \R^n$: The standard inner product on $\R^n$.

\item
$\left<X,Y\right>_{{\cal V}}$, $X,Y\in {\cal V}$: The inner product on ${\cal V}$.

\item
$\|X\|$, $X\in \R^n$: The standard norm on $\R^n$.

\item
$\|X\|_{{\cal V}}$: A norm on dimension-free vector space, or operator norm of operators over dimension-free vector space.

\item $\ltimes$: The left semi-tensor product (or M-product) of matrices.

\item $\lvtimes$: The left  vector product (or V-product ) of matrix with vector.

\item $\lplus$ ($\lminus$): The left M-addition (M-subtraction) of matrices.

\item $\lvplus$ ($\lvminus$): The left V-addition (V-subtraction) of vectors.

\item The set of all matrices:
$$
{\cal M}=\bigcup_{i=1}^{\infty}\bigcup_{j=1}^{\infty}{\cal M}_{i\times j}.
$$

\item The set of matrices:
$$
{\cal M}_{\mu}:=\left\{A\in {M}_{m\times n}\;\big|\; m/n=\mu \right\}.
$$
\item
$$
{\cal M}^{\mu}:=\oplus_{n=0}^{\infty}{\cal M}_{\mu^n}.
$$
\item $A\sim B$: Matrix equivalence (M-equivalence).

\item $\A$: M-equivalence class of $A$.

\item Set of vectors:
$$
{\cal V}=\bigcup_{n=1}^{\infty}{\cal V}_n.
$$

\item $x\lra y$: Vector equivalence (V-equivalence).

\item $\bar{x}$: V-equivalence class of $x$.

\item Quotient Matrix Space:
$$
\Sigma_{{\cal M}}={\cal M}/\sim.
$$
$$
\Sigma_{\mu}={\cal M}_{\mu}/\sim.
$$

\item Quotient Vector Space:
$$
\Omega={\cal V}/\lra.
$$

\item Set of Formal Polynomial on ${\cal M}$: ${\cal P}_{{\cal M}}$.

\item Set of Formal Polynomial on ${\cal M}^{\mu}$: ${\cal P}_{\mu}$.

\item Set of Formal Polynomial on $\Sigma$: ${\cal P}_{\Sigma}$.

\end{enumerate}

\vskip 2mm

\section{Cross-Dimensional S-Systems}

\subsection{Monoid of Matrices ${\cal M}$}

The set of all matrices is denoted as
\begin{align}\label{2.1.0}
{\cal M}=\bigcup_{m=1}^{\infty}\bigcup_{n=1}^{\infty}{\cal M}_{m\times n}.
\end{align}
Then ${\cal M}$ is the set of all (real) numbers, finite dimensional vectors, and matrices of arbitrary finite dimensions.

\begin{dfn}\label{d2.1.1} \cite{che11} Let $A,~B\in {\cal M}$. Precisely speaking, $A\in {\cal M}_{m\times n}$ and $B\in {\cal M}_{p\times q}$, and  $n\vee p=t$. Then the STP (or M-product) of $A$ and $B$ is defined as
\begin{align}\label{2.1.1}
A\ltimes B=\left(A\otimes I_{t/n}\right)\left(B\otimes I_{t/p}\right).
\end{align}
\end{dfn}

 The follow properties are fundamental \cite{che11,che12}:

\begin{prp}\label{p2.1.2}
\begin{enumerate}
\item (Associative Law)
\begin{align}\label{2.1.2}
(A\ltimes B)\ltimes C=A\ltimes (B\ltimes C).
\end{align}

\item (Distributive Law)
\begin{align}\label{2.1.3}
\begin{array}{l}
(\a A + \b B)\ltimes C=\a A\ltimes C+\b B\ltimes C;\\
C\ltimes (\a A+\b B)=\a C\ltimes A+\b C\ltimes B,\quad \a,\b\in \R.
\end{array}
\end{align}

\item
\begin{align}\label{2.1.4}
(A\ltimes B)^T=B^T\ltimes A^T.
\end{align}
\item Assume both $A$ and $B$ are invertible. Then
\begin{align}\label{2.1.5}
(A\ltimes B)^{-1}=B^{-1}\ltimes A^{-1}.
\end{align}

\end{enumerate}
\end{prp}

\begin{prp}\label{p2.1.3} $({\cal M},\ltimes)$ is a monoid.
\end{prp}

\noindent{\it Proof}.
Recall Definition \ref{d1.1}. The requirement (\ref{1.1}) comes from (\ref{2.1.2}). The identity is $1\in \R$.~\hfill$\Box$

\subsection{S-system on ${\cal V}$}

Let ${\cal V}_n$ be an $n$-dimensional real vector space. It can be briefly identified with $\R^n$ using the coordinates over an orthonormal basis.
Define
$$
{\cal V}:=\bigcup_{n=1}^{\infty} {\cal V}_n.
$$

\begin{dfn}\label{d2.2.1} Let $A\in {\cal M}_{m\times n}$ and $x\in {\cal V}_r$, and  $t=n\vee r$. Then the vector product (V-product) of $A$ with $x$ is defined as
\begin{align}\label{2.2.1}
A\lvtimes x:=\left(A\otimes I_{t/n}\right) \left(x\otimes \J_{t/r}\right).
\end{align}
\end{dfn}

\begin{rem}\label{r2.2.2}
\begin{enumerate}
\item V-product is the product of a matrix with a vector, STP is the product of a matrix with another matrix. For statement symmetry, STP is also called the M-product.
\item The V-product of a matrix $A$ with a vector $x$ can easily be extended to two matrices $A$ and $B$: Let $A\in {\cal M}_{m\times n}$, $B\in {\cal M}_{r\times s}$ and  $t=n\vee r$. Then
\begin{align}\label{2.2.101}
A\lvtimes B:=\left(A\otimes I_{t/n}\right) \left(B\otimes \J_{t/r}\right).
\end{align}
In this product the $B$ is considered as a set of $s$ column vectors in ${\cal V}_r$.
\item In fact a matrix $A\in {\cal M}_{m\times n}$ has two distinct meanings as follows: (i) It can be considered as a linear mapping from $\R^n$ to $\R^m$. This is the ``linear operator role" of a matrix. (ii) It can be considered as a subspace of $\R^m$, spanned by the columns of $A$. This is the ``subspace role" of a matrix. Then when we consider $A\times B$, it could be considered either a composed linear mapping of $A\circ B$ or the action of linear mapping $A$ on the subspace of $B$, precisely, the subspace spanned by the columns of $B$.  Fortunately, for standard case, the same conventional matrix product is applicable for these two functions. But now for generalized (dimension-free) matrix theory, these two functions need to be expressed by two different products. That is why we need M-product and V-product. They play the two different functions respectively.
\item From above argument one sees easily that both M-product and V-product are generalizations of conventional matrix product. That is, when $A$ and $B$ satisfy the dimension matching condition for conventional matrix product, both M- and V- products coincide with conventional one.
\end{enumerate}
\end{rem}

\begin{thm} \label{t2.2.3} Let the action of monoid $\left({\cal M},\ltimes\right)$ on ${\cal V}$ be defined by (\ref{2.2.1}). Then $({\cal M}, \lvtimes, {\cal V})$ is an S-system.
\end{thm}
\noindent{\it Proof.} The only thing, which is nontrivial, is to prove the following:
\begin{align}\label{2.2.3}
A\lvtimes (B\lvtimes x)=(A\ltimes B)\lvtimes x,\quad A,B\in {\cal M},\;x\in {\cal V}.
\end{align}
(\ref{2.2.3}) was proved in \cite{chepr}. \hfill $\Box$

Roughly speaking. the rest of this paper is devoted to this S-system and its variants. More topological and vector space structures will be posed on the monoid ${\cal M}$ and the objective state space ${\cal V}$. Then more dynamic properties of this S-system and its variants are revealed. Particularly, the corresponding S-system on quotient space is proposed and investigated.

\section{Linear Structure on  $\left({\cal M},\lvtimes, {\cal V}\right)$}

Since the S-system $\left({\cal M},\lvtimes, {\cal V}\right)$ is essentially the action of matrices on vector space, one of our purposes is to provide certain linear structures on the components of $\left({\cal M},\lvtimes, {\cal V}\right)$, such that this S-system becomes a ``linear system". This section will build these linear structures step by step.

\subsection{Vector Space Structure on ${\cal V}$}

To pose a vector space structure on ${\cal V}$, we need an addition on it. This has been done in \cite{chepr}.

\begin{dfn}\label{d3.1.1} \cite{chepr}
\begin{enumerate}
\item Let $x,~y\in {\cal V}$, precisely speaking, $x\in {\cal V}_m$, $y\in {\cal V}_n$ and $t=m\vee n$. Then the V-addition of $x$ and $y$, denoted by $x\lvplus y$, is defined as
\begin{align}\label{3.1.1}
x\lvplus y:=\left(x\otimes \J_{t/m}\right)+\left(y\otimes \J_{t/n}\right).
\end{align}
\item The corresponding V-substraction is defined as
\begin{align}\label{3.1.2}
x\lvminus y:=x\lvplus (-y).
\end{align}
\end{enumerate}
\end{dfn}

Now a natural question is: Is ${\cal V}$ with addition $\lvplus$ and conventional scalar product a vector space? To answer this question, we recall the definition of vector space first.

\begin{dfn}\label{d3.1.2} \cite{rad89} A set $V$ with an addition $+: V\times V\ra V$ and a scalar product $\cdot: \R\times V\ra V$ is called a vector space if the following are satisfied (for $x,y,z\in V$, $a,b\in \R$):
\begin{itemize}
\item[(I)] (i) $x+y\in V$; (ii) $x+y=y+x$; (iii) $(x+y)+z=x+(y+z)$; (iv) there exists a unique $\vec{0}$ such that $x+\vec{0}=x$; (v) there exists $-x$ such that $x+(-x)=\vec{0}$.
\item[(II)] (i) $ax\in V$; (ii) $a(bx)=(ab)x$; (iii) $(a+b)x=ax+bx$; (iv) $a(x+y)=ax+ay$; (v) $1\cdot x=x$; (vi) $0\cdot x=\vec{0}$; (vii) $a\cdot \vec{0}=\vec{0}$.
\end{itemize}
\end{dfn}

\begin{dfn}\label{d3.1.201} \cite{abr78} A set $V$ with an addition $+: V\times V\ra V$ and a scalar product $\cdot: \R\times V\ra V$ is called a pseudo-vector space if all the requirements for a vector space are satisfied except that the $\vec{0}$ is not unique.
\end{dfn}

In a pseudo-vector space, since $\vec{0}$ is not unique, it is clear that the ``inverse" $-x$ of $x$ may also not unique.
Then we need to make the concepts of $\vec{0}$ and $-x$ clear. Taking ${\cal V}$ as an example, we may define $\vec{0}$ being a set as
\begin{align}\label{3.1.3}
\vec{0}=\{{\bf 0}_r\;|\;r=1,2,\cdots\},
\end{align}
where ${\bf 0}_r$ is the zero of ${\cal V}_r$.

To deal with the problem on inverse $-x$ caused by multi-zero,  two ways are proposed:
\begin{itemize}
\item Assume $x\in {\cal V}_n$, then $-x\in {\cal V}_n$ and
\begin{align}\label{3.1.4}
-x + x={\bf 0}_n.
\end{align}
\item $-x\in {\cal V}$ and
\begin{align}\label{3.1.5}
-x \lvplus x\in \vec{0}.
\end{align}
\end{itemize}

Both (\ref{3.1.4}) and (\ref{3.1.5}) can not make ${\cal V}$ a vector space.

\noindent{\bf Case 1}:  We consider (\ref{3.1.4}) first. It leads to a new ``vector space" structure, called the hybrid vector space, which is defined as follows:

\begin{dfn}\label{d3.1.3} Let ${\cal S}_1$, $\cdots$,  ${\cal S}_n$ ($n$ is allowed to be $\infty$) be a set of vector spaces over $\R$.
$$
{\cal S}:=\bigcup_{i=1}^n {\cal S}_i
$$
with an addition $+_S$ and a scalar product $\cdot$ is called a hybrid vector space, if it satisfied the following conditions:
\begin{enumerate}
\item The overall addition $+_S$ satisfies the following:
\begin{itemize}
\item[(i)] I(i)-I(iii) of Definition \ref{d3.1.2};

\item[(ii)] when restricted $+_S$ to ${\cal S}_r$, it becomes $+_r$, which is the addition on $S_r$.
\end{itemize}

\item There is a set of zeros, denoted as $\vec{0}$, where
$$
\vec{0}\bigcap S_r ={\bf 0}_r,
$$
where ${\bf 0}_r$ is the zero of $S_r$.

\item If $x\in {\cal S}_r$, then $-x\in {\cal S}_r$ and $x+(-x)={\bf 0}_r$.

\item The scalar product $\cdot$, as a common scalar product of all component vector spaces, satisfies the following:

\begin{itemize}
\item[(i)] II(i)-II(iii), II(v) of Definition \ref{d3.1.2};
\item[(ii)] $a(x +_S y)=ax +_S ay$;
\item[(iii)] $0\cdot x\in {\bf 0}_r$,~for $x\in S_r$;
\item[(iv)] $a\cdot {\bf 0}_r={\bf 0}_r$.
\end{itemize}
\end{enumerate}
\end{dfn}

It is easy to verify the following result:

\begin{prp}\label{p3.1.4}  $({\cal V},\lvplus, \cdot)$ (where $\cdot$ is the conventional scalar product) with the inverse defined by (\ref{3.1.4}) is a hybrid vector space.
\end{prp}

\noindent{\bf Case 2}: Assume the inverse of any $x\in {\cal V}$ is defined by (\ref{3.1.5}).

Because of the non-uniqueness of zero, if we want all the other requirements of Definition \ref{d3.1.2} to be satisfied. Then one sees immediately that for each $x$ its inverse $-x$ is not unique. That is,
\begin{align}\label{3.1.6}
-x=\{y\;|\; x\lvplus y\in \vec{0} \}.
\end{align}
Eventually, we need
\begin{align}\label{3.1.7}
x=-(-x)=\{-y\;|\; y\in \{-x\}\}.
\end{align}
That is, each element has an equivalent set. That leads to a quotient space, which will be discussed latter.

\subsection{Topological Structure on Hybrid Vector Space ${\cal V}$}

In this subsection we propose a topological structure on the objective space ${\cal V}$ via a distance. First, we introduce an inner product on it.

\begin{dfn}\label{d3.2.1} \cite{mcd05}
Let $X$ be a vector space over $\R$. An inner product on $X$ is a function $\left<,\right>:X\times X\ra \R$ satisfying the following conditions for all $x,y,z\in X$ and $a,b\in \R$:
\begin{enumerate}
\item[(i)]
$$
\left<ax+by,z\right>=a \left<x,z\right>+b\left<y,z\right>.
$$
\item[(ii)]
$$
\left<x,y\right>=\left<y,x\right>.
$$
\item[(iii)] $\left<x,x\right>\geq 0$.
\item[(iv)] $\left<x,x\right>=0$ if and only if $x=0$.
\end{enumerate}
\end{dfn}

Consider the hybrid vector space ${\cal V}$. Can we define an inner product on it? The answer is ``Yes". We give the following definition:

\begin{dfn}\label{d3.2.1}
Let $x,~y\in {\cal V}$, precisely, $x\in {\cal V}_m$, $y\in {\cal V}_n$, and $t=m\vee n$. Then an inner product ${\cal V}\times {\cal V}\ra \R$ is defined as
\begin{align}\label{3.2.1}
\left<x,y\right>_{{\cal V}}:=\frac{1}{t}\left<(x\otimes \J_{t/m}),(y\otimes \J_{t/n})\right>.
\end{align}
\end{dfn}

\begin{rem}\label{r3.2.2}
\begin{enumerate}
\item It is easy to verify that (\ref{3.2.1}) is an inner product on ${\cal V}$. That is, the requirements in Definition \ref{d3.2.1} are all satisfied.
\item Since ${\cal V}$ is not exactly a vector space, the properties of inner product on a vector space should be verified carefully to see if each one is also valid for hybrid vector space.
\end{enumerate}
\end{rem}

Next, we define a norm on ${\cal V}$.
\begin{dfn}\label{d3.2.3}
Let $x\in {\cal V}_m\subset {\cal V}$. Then the norm of $x$ can be defined in a natural way as
\begin{align}\label{3.2.2}
\|x\|_{{\cal V}}:=\sqrt{\left<x,x\right>_{{\cal V}}}=\sqrt{\frac{1}{m}\left<x,x\right>}.
\end{align}
\end{dfn}

A straightforward verification shows that this norm is properly defined. That is, (i) $\|x\|\geq 0$ and $\|x\|=0$, if and only if $x=0$; (ii) $\|ax\|=|a|\|x\|$; (iii) $\|x+y\|\leq \|x\|+\|y\|$. (To prove the triangular inequality, we can first prove the Cauchy's inequality $|\left<x,y\right>|^2\leq \left<x,x\right>\left<y,y\right>$.  Refer to \cite{mcd05}.)

Using this norm, we can define the operator norm of a matrix  $A\in {\cal M}$ when it is acting on ${\cal V}$.

\begin{dfn}\label{d3.2.4}
Let $A\in {\cal M}_{m\times n}\subset {\cal M}$. Then the operator norm of $A$ is defined  as
\begin{align}\label{3.2.3}
\|A\|_{{\cal V}}:=\sup_{0\neq x \in {\cal V}} \frac{\|A\lvtimes x\|_{{\cal V}}}{\|x\|_{{\cal V}}}.
\end{align}
\end{dfn}

\begin{prp}\label{p3.2.5} Let $A\in {\cal M}_{m\times n}\subset {\cal M}$. Then
\begin{align}\label{3.2.4}
\|A\|_{{\cal V}}=\sqrt{\frac{n}{m}}\|A\|=\sqrt{\frac{n}{m}}\sqrt{\sigma_{max}(A^TA)}.
\end{align}
\end{prp}

\noindent{\it Proof}. First,
$$
\begin{array}{ccl}
\|A\|_{{\cal V}}&\geq &\|A\|_{{\cal V}_n}\\
~&=&\sup_{0\neq x \in {\cal V}_n} \frac{\sqrt{n}\|Ax\|}{\sqrt{m}\|x\|}\\
~&=&\sqrt{\frac{n}{m}}\sqrt{\sigma_{\max}(A^TA)}\\
~&=&\sqrt{\frac{n}{m}}\|A\|,
\end{array}
$$
where $\sigma(A)$ is the spectrum of $A$, $\|A\|$ is the standard norm of $A$.
(We refer to \cite{hor85} for last two equalities.)

On the other hand, let $x\in {\cal V}_r$. Then
$$
\begin{array}{ccl}
\|A\lvtimes x\|_{{\cal V}}&=&\sup_{x\in {\cal V}_r}\frac{\|(A\otimes I_{t/n})(x\otimes \J_{t/r})\|_{{\cal V}_t}}{\|x\otimes \J_{t/r}\|_{{\cal V}_t}}\\
~&\leq& \sup_{x\in {\cal V}_t}\frac{\|(A\otimes I_{t/n})x\|_{{\cal V}_t}}{\|x\|_{{\cal V}_t}}\\
~&=& \sqrt{\frac{t}{mt/n}}\|A\otimes I_{t/n}\|\\
~&=&\sqrt{\frac{n}{m}}\|A\|.
\end{array}
$$
The conclusion follows.
\hfill $\Box$

Using Definition \ref{d3.2.3} and Proposition \ref{p3.2.5}, a straightforward computation leads to the following result.

\begin{cor}\label{c3.2.501}
\begin{enumerate}
\item For a $\J_s$,  $s\in \N$, we have
\begin{align}\label{3.2.401}
\|x\otimes \J_s\|_{{\cal V}}=\|x\|_{{\cal V}}.
\end{align}
\item For an $I_s$,  $s\in \N$, we have
\begin{align}\label{3.2.402}
\|A\otimes I_s\|_{{\cal V}}=\|A\|_{{\cal V}}.
\end{align}
\end{enumerate}
\end{cor}

Finally, we can define a distance on ${\cal V}$ using the norm.

\begin{dfn}\label{d3.2.6} Let $x,~y\in {\cal V}$. Then the distance between $x$ and $y$ is defined as follows:
\begin{align}\label{3.2.5}
d(x,y):=\|x\lvminus y\|_{{\cal V}}.
\end{align}
\end{dfn}

\begin{rem}\label{r3.2.7}
\begin{enumerate}
\item Using the distance defined in (\ref{3.2.5}), we can get a topological structure on ${\cal V}$, which turns ${\cal V}$ into a topological space. Precisely, define $O^x_r:=\{y\;|\;d(x,y)<r\}$ and set $O:=\{O^x_r\;|\;x\in {\cal V}, r>0\}$. Using $O$ as a topological basis, we have the required topology ${\cal T}$.
\item Hereafter, we consider this topology ${\cal T}$ as the default topology on ${\cal V}$ and consider ${\cal V}$ as the topological space $({\cal V},{\cal T})$.
\item It is interesting to ask: if $({\cal V},d)$ a metric space? A matric space requires \cite{mcd05}:
\begin{itemize}
\item[(i)] $d(x,y)\geq 0$, with equality if and only if, $x=y$;
\item[(ii)] $d(x,y)=d(y,x)$;
\item[(iii)] $d(x,z)\leq d(x,y)+d(y,z)$.
\end{itemize}
It is easy to verify that  $({\cal V},d)$ satisfies all the three conditions except that $d(x,y)= 0$ does not imply $x=y$.
\item It is easy to show that $d(x,y)=0$ if and only if, there exist $\J_{\a}$ and $\J_{\b}$ such that
\begin{align}\label{3.2.6}
x\otimes \J_{\a}=y\otimes \J_{\b}.
\end{align}
This equation corresponds to the {\bf Case 2} of subsection (III.A).
\item The fact that $d(x,y)= 0$ does not imply $x=y$ leads to the fact that the topological space  $({\cal V},{\cal T})$ is not Hausdorff.
\end{enumerate}
\end{rem}

Next, we prove the continuity of the action of $A:{\cal V}\ra {\cal V}$.

\begin{prp}\label{p3.2.8} Given $A\in {\cal M}$. Then the mapping $x\mapsto A\lvtimes x$ is continuous.
\end{prp}

\noindent{\it Proof}. According to the operator norm of $A$, i.e., using (\ref{3.2.3}), we have
\begin{align}\label{3.2.6}
\|A\lvtimes x\|\leq \|A\|_{{\cal V}}\|x\|_{{\cal V}}.
\end{align}
Note that (\ref{3.2.4}) shows that $\|A\|_{{\cal V}}$ is independent of the dimension of $x$. Then the continuity is obvious.
\hfill $\Box$

\begin{cor}\label{c3.2.9} The S-system $\left({\cal M},\lvtimes,{\cal V}\right)$ is a pseudo dynamic system.
\end{cor}

\subsection{Hybrid Vector Space ${\cal M}_{\mu}$}

To find the vector space structure of matrices we split ${\cal M}$ as
\begin{align}\label{3.3.1}
{\cal M}=\bigcup_{\mu\in \Q_+}{\cal M}_{\mu},
\end{align}
where $\Q_+$ is the set of irreducible positive rational numbers. It is clear that (\ref{3.3.1}) is a partition.

Assume $\mu=\mu_y/\mu_x$, where $\mu_y\wedge \mu_x=1$. Then ${\cal M}_{\mu}$ can be further splat as
\begin{align}\label{3.3.2}
{\cal M}_{\mu}=\bigcup_{n=1}^{\infty}{\cal M}_{n\mu_y\times n\mu_x}.
\end{align}

On each ${\cal M}_{\mu}$ we can define an addition, called the matrix-addition (M-addition), as follows:

\begin{dfn}\label{d3.3.1} Let $A\in {\cal M}_{m\times n}\subset {\cal M}_{\mu}$, $B\in {\cal M}_{p\times q}\subset {\cal M}_{\mu}$ (i.e., $m/n=p/q=\mu$), and $t=m\vee p$. The M-addition of $A$ and $B$ is defined as
\begin{align}\label{3.3.3}
A\lplus B:=\left(A\otimes I_{t/m}\right)+\left(B\otimes I_{t/p}\right).
\end{align}
Correspondingly, the M-subtract is defined as
\begin{align}\label{3.3.4}
A\lminus B:=A\lplus (-B).
\end{align}
\end{dfn}

Similar to vector case, we have zero as a set of elements, as
$$
\vec{0}:=\left\{{\bf 0}_{n\mu_y\times n\mu_x\;|\; n=1,2,\cdots}\right\}.
$$

Similar argument as for vector case yields the following result:

\begin{prp}\label{p3.3.2} $\left({\cal M}_{\mu},\lplus, \cdot\right)$ is a hybrid vector space.
\end{prp}

\subsection{Linear S-System}

We give the following definition for a linear mapping, which is classical one for vector spaces \cite{rad89}:

\begin{dfn}\label{d3.4.1}
 Let ${\cal L},~{\cal X}$ be two (hybrid) vector spaces. A mapping $\varphi: {\cal L}\times {\cal X}\ra {\cal X}$ is called a linear mapping if
\begin{align}\label{3.4.1}
\begin{array}{ccl}
\varphi(L,ax+by)&:=&L(ax+by)=aLx+bLy,\\
~&~& L\in {\cal L}, \; x,y\in X,\;a,b\in \R.
\end{array}
\end{align}
\end{dfn}

\begin{dfn}\label{d3.4.2} Consider an S-system  $(S,\varphi,X)$.
\begin{enumerate}
\item $(S,\varphi,X)$ is called a linear S-system, if both $S$ and $X$ are vector spaces and $\varphi:S\times X\ra X$ is a linear mapping.
\item $(S,\varphi,X)$ is called a hybrid linear S-system, if both $S$ and $X$ are hybrid vector spaces and $\varphi:S\times X\ra X$ is a linear mapping.
\end{enumerate}
\end{dfn}

\begin{dfn}\label{d3.4.3} Assume $(G,\varphi,X)$ is a (linear, hybrid linear) S-system (or pseudo dynamic system,  dynamic system) and $H<G$ is a sub-monoid of $G$, then it is obvious that  $(H,\varphi,X)$ is also a (linear, hybrid linear) S-system (or pseudo dynamic system,  dynamic system). It will be called the sub (linear, hybrid linear) S-system (or pseudo dynamic system, dynamic system) of $(G,\varphi,X)$.
\end{dfn}

Back to the S-system $({\cal M}, \lvtimes, {\cal V})$, we have the following result:

Define
\begin{align}\label{3.4.2}
{\cal M}^{\mu}:=\bigcup_{n=0}^{\infty}{\cal M}_{\mu^n}.
\end{align}

Then it is clear that ${\cal M}^{\mu}< {\cal M}$ is a sub-monoid. Hence, we have the following result.

\begin{cor}\label{c3.4.4} $\left({\cal M}^{\mu},\lvtimes, {\cal V}\right)$ is a pseudo dynamic system.
\end{cor}

Unfortunately, ${\cal M}^{\mu}$ is not a (hybrid) vector space so far. We need to provide a hybrid vector space structure on it.

\section{Formal Polynomial}

\subsection{Direct Sum of Vector Spaces}

\begin{dfn}\label{d4.1.1} Assume $V_1$ and $V_2$ are two (hybrid) vector spaces over $\R$. The direct sum of $V_1$ and $V_2$ is defined as
\begin{align}\label{4.1.1}
V_1\oplus V_2:=x_1\oplus x_2,\quad x_1\in V_1,\;x_2\in V_2.
\end{align}
The addition ($+$) on $V_1 \oplus V_2$ is defined as
\begin{align}\label{4.1.2}
\begin{array}{l}
x_1\oplus x_2 + y_1\oplus y_2:=(x_1+_1y_1) \oplus (x_2+_2y_2),\\
\quad x_1,y_1\in V_1,\; x_2,y_2\in V_2,
\end{array}
\end{align}
where $+_i$ is the addition on $V_i$, $i=1,2$, respectively.

The scalar product is defined as
\begin{align}\label{4.1.3}
\begin{array}{l}
\lambda (x_1 \oplus x_2):=\lambda x_1 \oplus \lambda x_2,\\
\quad x_1\in V_1,\; x_2\in V_2,\;\lambda\in \R.
\end{array}
\end{align}
\end{dfn}

Then it is easy to prove the following result:

\begin{prp}\label{p4.1.2}
\begin{enumerate}
\item The direct sum of $n$ vector spaces is a vector space;\\
\item The direct sum of $n$ hybrid vector spaces is a hybrid vector space, where $n=\infty$ is allowed.
\end{enumerate}
\end{prp}

\begin{rem}\label{r4.1.201} The second property in Proposition \ref{p4.1.2} may also be considered as a definition. Consider
$$
{\cal M}=\oplus_{\mu\in \Q_+}{\cal M}_{\mu},
$$
the overall addition $\oplus$ is a formal operator. If we restrict it to each ${\cal M}_{\mu}$, it becomes $\lplus$. Since ${\cal M}_{\mu}$ is a hybrid vector space, if we restrict it to each vector space ${\cal M}_{m\times m\mu}$, it becomes the real vector space addition $+$. This restriction is slightly different from the original definition of hybrid vector space. So the hybrid vector space structure of a  direct sum of $n$ hybrid vector spaces should be understood in this way. That is, we have a set of hierarchy restrictions for the overall addition $\oplus$.
\end{rem}

\begin{exa}\label{e4.1.3}
Consider $W:={\cal M}_1\oplus {\cal M}_2\oplus {\cal V}$. According to Proposition \ref{p4.1.2}, it is a hybrid vector space.
Let $A_i\in {\cal M}_1$, $B_i\in {\cal M}_2$, $x_i\in {\cal V}$, $i=1,2$, be as follows:
$$
A_1=\begin{bmatrix}1&0\\0&-1\end{bmatrix};\quad
A_2=\begin{bmatrix}1&-1&2&-2\\
0&-1&1&0\\2&2&-1&1\\1&-1&0&-2\\
\end{bmatrix};
$$
$$
B_1=\begin{bmatrix}1&0\\0&2\\-2&1\\0&1\end{bmatrix};
\quad
B_2=
\begin{bmatrix}1\\2\end{bmatrix};
$$
$$
x_1=[1,0,-1]^T;\quad x_2=[-1,1]^T.
$$
Then
\begin{itemize}
\item[(i)]
$$
w_i=A_i\oplus B_i\oplus x_i\in W,\quad i=1,2.
$$
\item[(ii)]
$$
-w_i=-A_i\oplus -B_i\oplus -x_i.
$$
\item[(iii)]
$$
\begin{array}{ccl}
aw_1+bw_2&=&(aA_1\lplus bA_2)\oplus (aB_1\lplus bB_2)\\
~&~&\oplus (ax_1\lvplus bx_2)\\
~&=&\begin{bmatrix}a+b&-b&2b&-2b\\
0&a-b&b&0\\2b&2b&-a-b&b\\b&-b&0&-a-2b\\
\end{bmatrix}\\
~&~&\oplus \begin{bmatrix}a+b&0\\
0&2a+b\\2b-2a&a\\0&a+2b\\
\end{bmatrix}\oplus \begin{bmatrix}a-b\\
a-b\\-b\\b\\b-a\\b-a\\
\end{bmatrix}.
\end{array}
$$
\end{itemize}
\end{exa}

\begin{rem}\label{r4.1.4} A vector space can be considered as a special hybrid vector space which has only one component vector space. For statement ease hereafter  a hybrid vector space may either be vector space or a real hybrid vector space, which have more that one component vector spaces.
\end{rem}

\subsection{Formal Polynomials of Matrices}
\begin{dfn}\label{d4.2.1} Assume $W=\oplus_{\lambda\in \Lambda}S_{\lambda}$ is a direct sum of hybrid vector spaces $S_{\lambda}$, $\lambda\in \Lambda$. $W$ is called a formal polynomial, if the following are satisfied:
\begin{enumerate}
\item $\Lambda \subset \R_+$ is closed under product, that is, if $\lambda_i\in \Lambda$, $i=1,2$, then $\lambda_1\lambda_2\in \Lambda$.
\item  There exists a product $\odot:~W\times W\ra W$ satisfying ($w_1,~w_2,~w_3\in W$, $a,~b\in \R$)
\begin{itemize}
\item[(i)]
$$
\odot|_{S_{\lambda}}=\times_{\lambda}.
$$
\item[(ii)] (Linearity)
$$
\begin{array}{l}
(aw_1\oplus bw_2)\odot w_3=a(w_1\odot w_3)\oplus b(w_2\odot w_3)\\
w_1\odot (aw_2\oplus bw_3)=a(w_1\odot w_2)\oplus b(w_1\odot w_3).\\
\end{array}
$$
\item[(iii)]
$$
S_{\lambda_1}\odot S_{\lambda_2}\subset S_{\lambda_1\lambda_2}.
$$
\end{itemize}
\end{enumerate}
\end{dfn}

Define the formal polynomial of matrices as
\begin{align}\label{4.2.1}
{\cal P}_{{\cal M}}:=\oplus_{\mu\in \Q_+}{\cal M}_{\mu}.
\end{align}

Following facts are obvious:

\begin{prp}\label{p4.2.2}
\begin{enumerate}
\item ${\cal P}_{{\cal M}}$ is a hybrid vector space.
\item ${\cal M}\subset {\cal P}_{{\cal M}}$ is the set of monomials.
\item ${\cal P}_{{\cal M}}$ with product $\odot:=\ltimes$ is a set of formal polynomials.
\end{enumerate}
\end{prp}

Let $p\in {\cal P}_{{\cal M}}$. Then $p$ is called a formal polynomial of matrices.
Say,
$$
p=\oplus_{\mu\in \Q_+}A_{\mu},\quad A_{\mu}\in {\cal M}_{\mu}.
$$
Then we simply denote it as
\begin{align}\label{4.2.2}
p=\dsum_{\mu\in \Q_+}A_{\mu}z_{\mu},
\end{align}
where $z_{\mu}$ is a dummy variable satisfying formally $z_{\mu_1}z_{\mu_2}=z_{\mu_1\mu_2}$.

Following formulas are the consequence of the definitions of hybrid vector space and formal polynomial of matrices.

\begin{prp}\label{p4.2.3}
\begin{enumerate}
\item Let  $p,~q\in {\cal P}_{{\cal M}}$ with $p=\dsum_{\mu\in \Q^+}A_{\mu}z_{\mu}$ and $q=\dsum_{\mu\in \Q^+}B_{\mu}z_{\mu}$. Then the addition of $p$ and $q$ is
\begin{align}\label{4.2.3}
p\lplus q:=\dsum_{\mu\in Q_+}\left[A_{\mu} \lplus B_{\mu}\right]z_{\mu}.
\end{align}
The corresponding substraction is
\begin{align}\label{4.2.4}
p\lminus q:=\dsum_{\mu\in Q_+}\left[A_{\mu} \lplus (-B_{\mu})\right]z_{\mu}.
\end{align}

\item Let $p=\dsum_{\xi\in \Q^+}A_{\xi}z_{\xi}$, $q=\dsum_{\eta\in \Q^+}B_{\eta}z_{\eta}$. Then the product of $p$ and $q$ is
\begin{align}\label{4.2.5}
p\ltimes q=\dsum_{\mu\in Q_+}\left[\dsum_{\xi \eta =\mu}A_{\xi} \ltimes B_{\eta}\right]z_{\mu}.
\end{align}
\end{enumerate}
\end{prp}

We give an example to demonstrate  the computations on ${\cal P}_{\Sigma}$.

\begin{exa}\label{e4.2.4} Let
$$
p=A z_{1}+B z_{2};\quad q=C z_{1/2}+Dz_{1},
$$
where
$$
\begin{array}{ll}
A=\begin{bmatrix}
1&0\\1&1
\end{bmatrix};&
B=\begin{bmatrix}
1\\-1
\end{bmatrix};\\
C=\begin{bmatrix}
1&-1&2&-1\\
0&1&0&-1
\end{bmatrix};&
D=\begin{bmatrix}
1&-1&0\\0&2&1\\2&1&1
\end{bmatrix}.
\end{array}
$$

\begin{itemize}
\item[(1)]  Calculate
$p\lplus q$?

$$
p+q=A_1z_1+A_2z_2+A_{1/2}z_{1/2},
$$
where
$$
\begin{array}{ccl}
A_1&=&A\lplus D\\
~&=&\left(A\otimes I_3\right)+\left(D\otimes I_2\right)\\
~&=&\begin{bmatrix}
2&0&-1&0&0&0\\
0&2&0&-1&0&0\\
0&0&3&0&1&0\\
1&0&0&3&0&1\\
2&1&1&0&2&0\\
0&2&1&1&0&2
\end{bmatrix},
\end{array}
$$
$$
A_2=B,\quad A_{1/2}=C.
$$
\item[(2)]  Calculate
$p\ltimes q$?

$$
\begin{array}{l}
p\ltimes q\\
~=(A\ltimes C)z_{1/2}+(A\times D\lplus B\ltimes C)z_1+(B\ltimes D)z_2\\
~=\begin{bmatrix}1&-1&2&-1\\1&0&2&-2\end{bmatrix}z_{1/2}\\
~+\left[
\begin{array}{llllllll}
2&0&0&-1&-1&0&2&0\\
0&2&0&0&-1&-1&0&2\\
0&0&2&0&0&-1&-1&0\\
0&0&0&2&0&0&0&-1\\
0&0&0&0&3&0&0&0\\
0&0&0&0&0&3&0&0\\
0&0&0&1&-1&0&0&0\\
0&0&0&0&1&-1&0&0\\
2&0&0&0&1&1&-1&0\\
0&2&-1&1&0&1&0&-1\\
0&0&2&-1&2&0&1&0\\
0&0&0&2&-1&2&0&1\\
\end{array}\right.\\
~~~~\left.\begin{array}{llll}
0&-1&0&0\\
0&0&-1&0\\
2&0&0&-1\\
0&-1&0&0\\
1&0&-1&0\\
0&1&0&-1\\
0&2&1&0\\
0&0&2&1\\
-1&0&0&2\\
0&2&0&0\\
1&0&2&0\\
0&1&0&2\\
\end{array}\right]z_1
+\begin{bmatrix}1&-1&0\\
0&2&1\\2&1&1\\-1&1&0\\0&-2&-1\\
-2&-1&-1\end{bmatrix}z_2.
\end{array}
$$
\end{itemize}
\end{exa}

Recall ${\cal M}^{\mu}$ defined by (\ref{3.4.2}). We can merge it into ${\cal P}_{{\cal M}}$ as a subspace of hybrid vector space. Define
\begin{align}\label{4.2.6}
{\cal P}_{\mu}:=\left\{ \dsum_{i=0}^{\infty}A_iz^i\;|\; A_i\in {\cal M}_{\mu^i}, z^i:=z_{\mu^i}\right\}.
\end{align}

Then it is easy to verify that ${\cal P}_{\mu}$ with M-product $\ltimes$ is a formal polynomial over a hybrid vector space.
$p\in {\cal P}_{\mu}$ is expressed exactly as a polynomial and the addition and product of $p,~q\in {\cal P}_{\mu}$ are also exactly the same as the ones for polynomials.

Particularly, if a formal polynomial in (\ref{4.2.6}) is generated by a single matrix $A$, that is, $A_i=c_i A^i$, $i=0,1,\cdots$, for a fixed $A$, then such kind of formal polynomials form a subset of ${\cal P}_{\mu}$, denoted as
\begin{align}\label{4.2.7}
{\cal P}(A):=\left\{\dsum_{i=0}^{\infty} c_i A^i z^{i}\;|\; c_i\in \R\right\}.
\end{align}
$p\in {\cal P}(A)$ is called a principle formal polynomial (PFP). Principle formal polynomials are of particular interest. We give the following examples to see that the analytic functions of matrices can be extended to non-square matrices.

\begin{exa}\label{e4.2.5}
\begin{enumerate}
\item Let $f(x)$ be an analytic real function with its Taylor expansion as
\begin{align}\label{4.2.8}
f(x)=\dsum_{n=0}^{\infty} c_nx^n.
\end{align}
Then for any $A\in {\cal M}$ we can define a principle formal polynomial as
\begin{align}\label{4.2.9}
f(A)=\dsum_{n=0}^{\infty} c_nA^n.
\end{align}
\item For any $A\in {\cal M}$ we have
\begin{align}\label{4.2.10}
e^A=\dsum_{n=0}^{\infty} \frac{1}{n!}A^nz^{n}\in {\cal P}(A).
\end{align}
\item Similarly we also have
\begin{align}\label{4.2.11}
\sin(A):=\dsum_{k=1}^{\infty}(-1)^{k-1}\frac{1}{(2k-1)!}A^{2k-1}z^{2k-1}\in {\cal P}(A).
\end{align}
\begin{align}\label{4.2.12}
\cos(A):=\dsum_{k=0}^{\infty}(-1)^{k}\frac{1}{(2k)!}A^{2k}z^{2k}\in {\cal P}(A).
\end{align}
\item A straightforward computation shows the Euler formula
\begin{align}\label{4.2.13}
e^{iA}=\cos(A)+i\sin(A),\quad A\in {\cal M}.
\end{align}
\end{enumerate}
\end{exa}

Later on, one will see that the formal polynomial of matrices plays a very important role for describing the trajectories of cross-dimensional linear systems.

\section{Discrete Time Linear System}

\subsection{Discrete-time Linear Weak Dynamic System on ${\cal V}$}

Consider a discrete-time S-system as
\begin{align}\label{5.1.1}
\begin{array}{l}
x(t+1)=A(t)\lvtimes x(t),\quad \mbox{with} ~x(0)=x_0,\\
\quad x(t)\in {\cal V},\; A(t)\in {\cal M}.
\end{array}
\end{align}

First, it is clear that this is a well posed evolutionary system. Because for a given $x_0$ the trajectory $\{x(t)\;|\;t=0,1,\cdots\}$ can be determined iteratively and uniquely. Second, the system is evolving at each step as an S-system $\left({\cal M},\lvtimes,{\cal V}\right)$, we well classify the system (\ref{5.1.1}) by the characteristics of $\left({\cal M},\lvtimes,{\cal V}\right)$. Then we have the following result.

\begin{prp}\label{p5.1.1} The system (\ref{5.1.1}) is a hybrid linear pseudo dynamic system.
\end{prp}

\noindent{\it Proof}.

It is well known that ${\cal V}$ is a hybrid vector space. Since
$$
{\cal M}=\oplus_{\mu\in \Q_+}{\cal M}_{\mu},
$$
according to Proposition \ref{p4.1.2} ${\cal M}$ is also a hybrid vector space. Corollary \ref{c3.2.9} shows that it is a pseudo dynamic system. Finally, it is easy to verify that it is linear according to Definition \ref{d3.4.2}.
\hfill $\Box$.

We give an example to depict the trajectory of system (\ref{5.1.1}).

\begin{exa}\label{e5.1.2}
Consider system (\ref{5.1.1}) with
$$
A(t)=\begin{bmatrix}
\sin(\frac{\pi t}{2})&0&1&-\cos(\pi t)\\
-1&\cos(\frac{\pi t}{2})&\sin(\frac{\pi (t+1)}{2})&1
\end{bmatrix},
$$
and $x_0=[1,0,-1]^T\in \R^3$. Then it is easy to calculate that
$$
A(0)=\begin{bmatrix}
0&0&1&-1\\
-1&1&1&1
\end{bmatrix};\quad
A(1)=\begin{bmatrix}
    1&0&1&1\\
    -1&0&0&1
    \end{bmatrix}
$$
$$
A(2)=\begin{bmatrix}
0&0&1&-1\\
-1&-1&-1&1
\end{bmatrix};\quad
A(3)=\begin{bmatrix}
-1&0&1&1\\
-1&-0&-0&1
\end{bmatrix}.
$$

In general,
$$
A(k)=A(i),\quad mod(k,4)=i.
$$
Then it is easy to calculate that
$$
x(1)=A(0)\lvtimes x_0=[1,1,0,-1,-2,-3]^T.
$$
Similarly, we can calculate that
$$
\begin{array}{l}
x(2)=[-2,-3,-4,-3,-4,-4]^T\\
x(3)=[1,1,0,4,5,7]^T\\
x(4)=[8,10,11,4,6,6]^T\\
x(5)=[-2,-2,0,12,13,13]^T\\
x(6)=[23,23,24,15,15,15]^T\\
x(7)=[0,0,0,-46,-47,-47]^T\\
x(8)=[-93,-93,-94,-47,-47,-47]^T\\
x(9)=[0,0,0,-94,-95,-95]^T\\
x(10)=[-189,-189,-190,-95,-95,-95]^T,\cdots
\end{array}
$$
In fact, we can figure out that
$$
B_i:=A(i)|_{\R^6}, i=0,1,\cdots.
$$
Then we have (We refer reader to the next subsection for details of the calculation.)
$$
B_0=\begin{bmatrix}
0&0&0&1&-1&0\\
0&0&0&1&0&-1\\
0&0&0&0&1&-1\\
-1&1&0&1&1&0\\
-1&0&1&1&0&1\\
0&-1&1&0&1&1
\end{bmatrix};
$$
$$
B_1=\begin{bmatrix}
1&0&0&1&1&0\\
1&0&0&1&0&1\\
0&1&0&0&1&1\\
-1&0&0&0&1&0\\
-1&0&0&0&0&1\\
0&-1&0&0&0&1\\
\end{bmatrix};
$$
$$
B_2=\begin{bmatrix}
0&0&0&1&-1&0\\
0&0&0&1&0&-1\\
0&0&0&0&1&-1\\
-1&-1&0&-1&1&0\\
-1&0&-1&-1&0&1\\
0&-1&-1&0&-1&1\\
\end{bmatrix};
$$
$$
B_3=\begin{bmatrix}
-1&0&0&1&1&0\\
-1&0&0&1&0&1\\
0&-1&0&0&1&1\\
-1&0&0&0&1&0\\
-1&0&0&0&0&1\\
0&-1&0&0&0&1\\
\end{bmatrix}.
$$
In general,
$$
B(k)=B_i,\quad mod(k,4)=i,\quad k=0,1,\cdots,
$$
and
$$
x(t+1)=B(t)x(t),\quad t\geq 1.
$$
\end{exa}

\subsection{Time Invariant Linear System}

\begin{dfn}\label{d5.2.1} Consider system (\ref{5.1.1}).
When $A(t)=A$, $t\geq 0$, we have
\begin{align}\label{5.2.1}
\begin{array}{l}
x(t+1)=A\lvtimes x(t),\quad \mbox{with} ~x(0)=x_0,\\
\quad x(t)\in {\cal V},\; A\in {\cal M},
\end{array}
\end{align}
which is called a time invariant discrete time linear pseudo dynamic system.
\end{dfn}

\begin{dfn}\label{d5.2.2}
Let $A\in {\cal M}$. ${\cal V}_r$ is $A$-invariant if
$$
A\lvtimes x\in {\cal V}_r,\quad \forall x\in {\cal V}_r.
$$
\end{dfn}

\begin{dfn}\label{d5.2.3}
$A$ is called a dimension-bounded operator, if for any $x(0)=x_0\in {\cal V}_r$, there exist an $N\geq 0$ and an $s$, such that $x(t)\in {\cal V}_{r_*}$, $t\geq N$.
\end{dfn}

Note that if at certain time $N\geq 0$ the trajectory reaches ${\cal V}_{r_*}$, that is, $x(N)\in {\cal V}_{r_*}$, moreover, if ${\cal V}_{r_*}$ is an $A$-invariant space, then we have $x(t)\in {\cal V}_{r_*}$, $t\geq N$. If the invariant spaces exist, does it necessary that each trajectory enters to an invariant space? Also, if both $x_0$ and $x_0'$ belong to the same  ${\cal V}_r$, will their trajectories both enter or not enter to the same $A$-invariant space? In fact, the answers to these two problems are ``Yes". It was proved that the  $A$-invariant spaces have very nice properties.

\begin{prp}[\cite{chepr}]\label{p5.2.4} $A$ is dimension-bounded, if and only if, $A\in {\cal M}_{\mu}$, where $\mu=\mu_y/\mu_x$, $\mu_y\wedge \mu_x=1$ and $\mu_y=1$.
\end{prp}

For instant $A\in {\cal M}_{3\times 6}$ is dimension-bounded, $A\in {\cal M}_{3\times 5}$ is not dimension-bounded.

Next, for an dimension-bounded $A\in {\cal M}_{k\times k\mu_x}$ and a starting point $x_0\in {\cal V}_{r_0}$, we try to find the dimensions of its trajectory and the invariant space, to which it will enter. First, we define a sequence of natural numbers $r_i$ recursively as
\begin{align}\label{5.2.2}
\begin{array}{l}
~~~\xi_0:=\frac{r_0\vee (k\mu_x)}{k\mu_x},\\
\begin{cases}
r_i=\xi_{i-1}k\\
\xi_i=\frac{\xi_{i-1}\vee \mu_x}{\mu_x},\quad i=1,2,\cdots.
\end{cases}
\end{array}
\end{align}
Then we can prove the dimension of $x(t)$ is equal to $r_t$.

\begin{prp}\label{p5.2.5} Assume $A$ is dimension-bounded and $x_0\in {\cal V}_{r_0}$. Then
\begin{enumerate}
\item the dimension of $x(t)$ is $r_t$, i.e., $x(t)\in {\cal V}_{r_t}$, where $r_t$ are calculated recursively by (\ref{5.2.2});
\item there exists a (smallest) $i_*\geq 0$, such that
\begin{align}\label{5.2.3}
r_{i}=r_{i_*}:=r_*,\quad i\geq i_*.
\end{align}
\end{enumerate}
\end{prp}

\noindent{\it Proof}.
\begin{enumerate}
\item According to (\ref{5.2.2}),
$$
r_0\vee k\mu_x=\xi_0 k\mu_x.
$$
Then
$$
\dim(x(1))=k\frac{r_0\vee k\mu_x}{k\mu_x}=k\xi_0:=r_1.
$$
Next, since
$$
r_1\vee k\mu_x=k(\xi_0\vee \mu_x)=k\xi_1\mu_x,
$$
then
$$
\dim(x(2)):=r_2=k\xi_1.
$$
Similarly, we have $\dim(x(k))=r_k$, $k\geq 1$.
\item By definition it is clear that
$$
\xi_{i+1}\leq \xi_i, \quad i\geq 1,
$$
and as long as $\xi_{i_*+1}=\xi_{i_*}$, all $\xi_i=\xi_{i_*}$, $i>i_*$. The existence of $i_*$ is obvious.
\end{enumerate}
\hfill $\Box$

\begin{cor}\label{c5.2.6} Assume $A$ is dimension-bounded, precisely, $A\in {\cal M}_{k\times k\mu_x}$. ${\cal V}_r$ is $A$ invariant, if and only if,
\begin{align}\label{5.2.4}
r=\frac{r\vee (k\mu_x)}{\mu_x}.
\end{align}
\end{cor}

\begin{cor}\label{c5.2.7} Assume $A$  is dimension-bounded. Then $A$  has countable many invariant spaces ${\cal V}_r$.
\end{cor}

\noindent{\it Proof}. The existence of invariant space comes from Proposition \ref{p5.2.5}. We have only to prove that the number of invariant spaces is infinite. It will be proved by contradiction.  Assume there are only finite of them, and the largest dimension of invariant spaces is $r_{\max}$. Let $s>1$ and $s\wedge k\mu_x=1$ (say, $s=k\mu_x+1$). Then
$$
\frac{k\mu_x\vee (sr_{\max})}{\mu_x}=\frac{s(k\mu_x\vee r_{\max})}{\mu_x}=sr_{\max}.
$$
According to Corollary \ref{c5.2.6}, ${\cal V}_{sr_{\max}}$ is $A$-invariant, which contradicts to the maximality of $r_{\max}$.
\hfill $\Box$

Assume ${\cal V}_{r}$ is $A$ invariant space. Then $A|_{{\cal V}_r}$ is a linear mapping on ${\cal V}_{r}$. Hence there exists a square matrix, denoted by $A_r$, such that
$$
A|_{{\cal V}_r}=A_r.
$$
It is easy to calculate $A_r$ as
\begin{align}\label{5.2.5}
\Col_i(A_r)=A\lvtimes \d_r^i,\quad i=1,2,\cdots,r.
\end{align}

\begin{rem}\label{r5.2.8} Assume $A\in {\cal M}_{k\times k\mu_x}$ is dimension-bounded. From Corollary \ref{c5.2.7} one
sees easily that there is no largest $A$ invariant space. But $A$ has smallest (non-zero) invariant space. It is easy to verify that the smallest invariant space of $A$ is ${\cal V}_k$.
Then we denote
\begin{align}\label{5.2.6}
A|_{{\cal V}_k}:=A_*.
\end{align}
\end{rem}

Finally, we provide a simple formula to calculate $A_r$.

\begin{prp}\label{p5.2.9} Assume $A\in {\cal M}_{k\times k\mu_x}$ and $r>0$ satisfies  (\ref{5.2.4}), then
\begin{align}\label{5.2.7}
A_r:=A|_{{\cal V}_r}=\left(A\otimes I_{r/k}\right)\left(I_r\otimes \J_{\mu_x}\right).
\end{align}
\end{prp}

\noindent{\it Proof}. According to (\ref{5.2.5}), one sees that
$$
A_r=A\lvtimes I_r.
$$
Using (\ref{5.2.4}), we have
$$
\begin{array}{ccl}
A\lvtimes I_r&=&\left(A\otimes I_{\frac{r\vee k\mu_x}{k\mu_x}}\right)\left(I_r\otimes \J_{\frac{r\vee k\mu_x}{r}}\right)\\
~&=&\left(A\otimes I_{\frac{r}{k}}\right)\left(I_r\otimes \J_{\mu_x}\right)\\
\end{array}
$$
\hfill $\Box$

\subsection{Solution of Linear System}

This subsection considers the calculation of the trajectory of (\ref{5.2.1}). We assume $A$ is dimension-bounded. In fact, it is easy if $A$ is not dimension-bounded, then stating from any $x_0$ the dimension $r(t)$ of $x(t)$ will go to $\infty$.

\begin{prp}\label{p5.3.1} Consider system (\ref{5.2.1}). Assume $x(0)=x_0\not\in \vec{0}$ and $x(t)\in {\cal V}_{r_t}$. Then
\begin{align}\label{5.3.1}
\lim_{t\ra \infty}r_t=\infty.
\end{align}
\end{prp}
\noindent{\it Proof}. Assume $A\in {\cal M}_{\mu}$, where $\mu=\mu_y/\mu_x$ and $\mu_y\wedge \mu_x=1$. According to Proposition \ref{p5.2.4} $\mu_y>1$. Now let $x_0\in {\cal V}_{r_0}$, then the dimension sequence $\{r_t\;|\;t=0,1,\cdots\}$ is independent of the special $x_0$. That is, stating from any $\xi_0\in {\cal V}_{r_0}$, the dimension sequence is the same.

Now assume $\{r_t\;|\;t=0,1,\cdots\}$ is bounded. Using Proposition \ref{p5.2.4} again, the sequence does not have fixed point. Hence,  there is at least one cycle, say, $(r_p,r_{p+1},\cdots,r_{p+\ell}=r_p)$. It follow that ${\cal V}_{r_p}$ is an $A^{\ell}$-invariant space. But $A^{\ell}\in {\cal M}_{\mu_y^{\ell}\times \mu_x^{\ell}}$ and $\mu_y^{\ell}\neq 1$. This is a contradiction. Then we have

\begin{align}\label{5.3.2}
\overline{\lim}_{t\ra \infty}r_t=\infty.
\end{align}

Note that
$$
r_{t+1}=m\frac{n\vee r_t}{n}.
$$
It follows that if $r_{t+1}\geq r_t$ then $r_{t+2}\geq r_{t+1}$. Hence as long as $r_{t_0+1}\geq r_{t_0}$, then $r_t,\;t\geq t_0$ is monotonically nondecreasing. Because of (\ref{5.3.2}), such $t_0$ exists. Hence, the limit for $t\ra \infty$ exists. Then (\ref{5.3.2}) implies (\ref{5.3.1}).
\hfill $\Box$

proposition \ref{p5.3.1} shows that it is reasonable to  assume $A$ being dimension-bounded.

Using invariant space, the trajectory of (\ref{5.2.1}), starting from $x_0\in {\cal V}$, can be calculated. It is described by the following example.

\begin{exa}\label{e5.3.2} Consider system (\ref{5.2.1}), where
$$
A=
\begin{bmatrix}
1&2&-1&0\\
1&-2&2&-1
\end{bmatrix}.
$$
\begin{enumerate}
\item If $\bar{x}(0)=\overline{x_0}$ and $x_0=[1,2,-1]$, find the solution:

Consider the minimum realization with $A$ and $x_0$ as in the above.

It is easy to calculate that $r_*=r_1=6$. Then we have
$$
x(1)=A\lvtimes x_0=[1,3,6,4,2,-4]^T.
$$
Moreover,
$$
A_*=A|_{{\R^6}}
=\begin{bmatrix}
     1&     2&     0&    -1&     0&     0\\
     1&     0&     2&    -1&     0&     0\\
     0&     1&     2&     0&    -1&     0\\
     1&    -2&     0&     2&    -1&     0\\
     1&     0&    -2&     2&     0&    -1\\
     0&     1&    -2&     0&     2&    -1
\end{bmatrix}.
$$
Then
$$
x(t)=A_*^{t-1}x(1),\quad t\geq 2.
$$
For instance, $x(2)=A_*x(1)=[3,9,13,1,1,-1]^T$, $x(3)=A_*x(2)=[20,28,34,-14,-20,-14]^*$, etc.

\item If $x_0=[1,0,2,-2,-1,1,2,0]$, find the solution:

It is easy to calculate that $r_1=4$, $r_2=r_*=2$. Then we have
$$
\begin{array}{l}
x(1)=A\lvtimes x_0=[6,-5,-7,6]^T,\\
x(2)=A\lvtimes x(1)=[3,-4]^T.
\end{array}
$$
$$
A_*=A|_{r_*}=\begin{bmatrix}
3&    -1\\
-1&    1
\end{bmatrix}.
$$
Hence
$$
x(t)=A_*^{t-2}x(2),\quad t\geq 3.
$$
For instance, $x(3)=A_*x(2)=[13,-7]^T$, $x(4)=A_*x(3)=[46,-20]^T$, etc.
\end{enumerate}
\end{exa}

\section{Continuous Time Linear System}

\subsection{Action of ${\cal P}_{{\cal M}}$ on ${\cal V}$}

Consider a formal polynomial
\begin{align}\label{6.1.1}
p=\dsum_{\mu\in \Q_+}A_{\mu}z_{\mu}\in {\cal P}_{{\cal M}}.
\end{align}
To make it a proper operator we assume the absolute convergence of formal polynomials:

\vskip 2mm

\noindent{\bf A-1}:
\begin{align}\label{6.1.2}
\dsum_{\mu\in \Q_+}\|A_{\mu}\|_{{\cal V}}<\infty.
\end{align}

\vskip 2mm

To see that the assumption {\bf A-1} is reasonable we consider
$$
e^A=\dsum_{n=0}^{\infty}\frac{1}{n!}A^n
$$
Then it is clear that
$$
\dsum_{n=0}^{\infty} \frac{1}{n!}\|A^n\|_{{\cal V}}\leq \dsum_{n=0}^{\infty} \frac{1}{n!}\|A\|^n_{{\cal V}}<\infty.
$$

Hereafter, {\bf A-1} is always assumed. Precisely speaking, hereafter, ${\cal P}_{{\cal M}}$ is used for the set of formal polynomials which satisfy {\bf A-1}.

\begin{dfn}\label{d6.1.1} Assume $p\in {\cal P}_{{\cal M}}$ is as in (\ref{6.1.1}). Then the action of  ${\cal P}_{{\cal M}}$ on ${\cal V}$ is defined as
\begin{align}\label{6.1.3}
p\lvtimes x:=\dsum_{\mu\in \Q_+}A_{\mu}\lvtimes x,\quad x\in {\cal V}.
\end{align}
\end{dfn}

Using this definition, it is easy to see the following:

\begin{prp}\label{p6.1.2} $\left( {\cal P}_{{\cal M}},\lvtimes, {\cal V}  \right)$ is a pseudo dynamic linear system.
\end{prp}

\subsection{Solution to Continuous Time Linear System}

\begin{dfn}\label{d6.2.1} A continuous-time linear system is defined as
\begin{align}\label{6.2.1}
\begin{cases}
\dot{x}(t)=A(t) \lvtimes x(t), \quad x(0)=x_0\\
\quad x(t)\in {\cal V}, A(t)\in {\cal M}.
\end{cases}
\end{align}
\end{dfn}

It is obvious that (\ref{6.2.1}) is a generalization of the classical linear system, which has square $A(t)$. As $A(t)=A$, the system is called a time-invariant linear system. Even if for classical linear system, if $A(t)$ is time-varying, in general, it is difficult to find a closed form solution. Hence we consider only the time-invariant case.

Using formal polynomial, we can easily obtain the solution of (\ref{6.2.1}) for $A(t)=A$:

\begin{prp}\label{p6.2.2} Assume $A(t)=A$, the solution of (\ref{6.2.1}) is
\begin{align}\label{6.2.2}
x(t)=e^{A t}\lvtimes x_0.
\end{align}
\end{prp}

\noindent{\it Proof}. By definition
$$
e^{At}\lvtimes x_0=\dsum_{n=0}^{\infty}\frac{t^n}{n!}[A^n\lvtimes x_0].
$$
A straightforward verification shows the conclusion.
\hfill $\Box$

Now the problem is how to calculate this trajectory.

{\bf Case 1} $A$ is not dimension-bounded. In this case  the solution will have unbounded dimension. Since {\bf A-1} is satisfied, we can find finite dimensional approximated solution up to any accuracy. This case may represent some diffusion process. It is worth being discussed in the further.

{\bf Case 2} $A$ is dimension-bounded. In this case the closed form solution can be obtained. In the following we consider how to calculate the trajectory.

Assume $x_0\in {\cal V}_{r_0}$. Since $A$ is dimension-bounded, using (\ref{5.2.2}) we can find $r_s=r_*$ such that ${\cal V}_{r_*}$ is $A$-invariant space. Then we can calculate the trajectory as

\begin{align}\label{6.2.3}
\begin{array}{ccl}
e^{At}\lvtimes x_0&=&\,x_0\lvplus tA\lvtimes x_0\lvplus\frac{t^2}{2!}A^2\lvtimes x_0\lvplus\cdots\lvplus\frac{t^s}{s!}A^s\lvtimes x_0\lvplus\cdots\\
                  &=&\,\left[x_0\lvplus tx_1\lvplus\frac{t^2}{2!}x_2\lvplus\cdots\lvplus\frac{t^{s-1}}{(s-1)!}x_{s-1}\right]\\
                   & &\,\lvplus\left[\frac{t^s}{s!}I_{r_*}+\frac{t^{s+1}}{(s+1)!}A_*+\cdots \right]x_s\\
                 &:=&\,I\lvplus II.
\end{array}
\end{align}
where $x_j=A^j\lvtimes x_0$, $j=1,2,\cdots,s$.

First, we assume $A_*$ is invertible. Then we can convert $II$ into $e^{A_*t}$ by adding finite terms. It results in
\begin{align}\label{6.2.4}
\begin{array}{l}
e^{At}\lvtimes x_0=\left[x_0\lvplus tx_1\lvplus\frac{t^2}{2!}x_2\lvplus\cdots\lvplus\frac{t^{s-1}}{(s-1)!}x_{s-1}\right]\\
~~\lvminus A_*^{-s}\left[I_{r_*}+tA_*+\frac{t^2}{2!}A_*^2+\cdots+\frac{t^{s-1}}{(s-1)!}A_*^{s-1}
\right]x_s\\
~~\lvplus A_*^{-s} e^{ A_*  t}x_s.\\
\end{array}
\end{align}

Second, we consider the general case. Note that
\begin{align}\label{6.2.5}
\begin{cases}
\frac{d^s}{dt^s}II=e^{A_*t}x_s\\
\frac{d^j}{dt^j}II\big|_{t=0}=0,\quad 0\leq j<s.
\end{cases}
\end{align}
We, therefore, have the solution of (\ref{6.2.5}) as
\begin{align}\label{6.2.6}
II=\int_0^td\tau_1\int_0^{\tau_1}d\tau_2\cdots\int_0^{\tau_{s-1}}e^{A_*\tau_s}d\tau_sx_s.
\end{align}
Note that using Jordan canonical form \cite{hor85}, (\ref{6.2.6}) can be calculated directly.

We use a numerical example to depict the calculation.

\begin{exa}\label{e6.2.3} Let
\begin{align}\label{6.2.7}
A=\begin{bmatrix}
1&2&-1&0\\
1&-2&2&-1
\end{bmatrix}.
\end{align}
Find $x(t)$ for $x_0$ given as follows.
\begin{enumerate}
\item
Assume $x_0=[1,2,-1]^T\in \R^3$:

Since $x_0\in \R^{3}$, $x_1=x(1)\in \R^6$, then $r_*=r_1=6$.
It is easy to calculate that
$$
x_1=[1,3,6,4,2,-4]^T;
$$
and
\begin{align}\label{6.2.8}
A_*=
\begin{bmatrix}
1&2&0&-1&0&0\\
1&0&2&-1&0&0\\
0&1&2&0&-1&0\\
1&-2&0&2&-1&0\\
1&0&-2&2&0&-1\\
0&1&-2&0&2&-1
\end{bmatrix};
\end{align}
$$
A_*^{-1}=~~~~~~~~~~~~~~~~~~~~~~~~~~~~~~~~~~~~~~~~~~~~~~~~~~~~~~~~~
$$
$$
\begin{bmatrix}
0.5&0.1667&-0.1667&0.5&-0.1667&0.1667\\
0.5&-0.5&0.5&0.5&-0.5&0.5\\
0&0&0.5&0.5&-0.5&0.5\\
0.5&-0.8333&0.8333&1.5&-1.1667&1.1667\\
0.5&-0.5&0.5&1.5&-1.5&1.5\\
1.5&-1.5&0.5&2.5&-2.5&1.5\\
\end{bmatrix}.
$$
Finally, we have
\begin{align}\label{6.2.9}
x(t)=x_0\lvminus A_*^{-1}x_1 \lvplus A_*^{-1}e^{A_*t}x_1.
\end{align}

\item
Assume $x_0=[1,0,2,-2,-1,1,2,0]^T\in \R^8$:

It is easy to calculate that
$$
\begin{aligned}
x_1&=A\lvtimes x_0=[6,-5,-7,6]^T\in \R^4;\\
x_2&=A\lvtimes x_1=[3,-4]^T\in\R^2;\\
r_*&=r_2;\\
\end{aligned}
$$
\begin{align}\label{6.2.10}
A_*=\begin{bmatrix}
3&-1\\
-1&1
\end{bmatrix},\quad
A_*^{-1}=\frac{1}{2}\begin{bmatrix}
1&1\\
1&3
\end{bmatrix}.
\end{align}

We use two methods to calculate $x(t)$:
\begin{itemize}
\item Using formula (\ref{6.2.4}), we have
\begin{align}\label{6.2.11}
x(t)=x_0\lvplus x_1t\lvplus [A_*^{-2}e^{A_*t}-A_*^{-2}-tA_*^{-1}]x_2,
\end{align}
where
$$
e^{A_*t}=\frac{1}{2\sqrt{2}}\begin{bmatrix}
a&b\\c&d
\end{bmatrix},
$$
with
$$
\begin{array}{l}
a=(\sqrt{2}+1)e^{(2+\sqrt{2})t}+(\sqrt{2}-1)e^{(2-\sqrt{2})t},\\
b=-e^{(2+\sqrt{2})t}+e^{(2-\sqrt{2})t},\\
c=-e^{(2+\sqrt{2})t}+e^{(2-\sqrt{2})t},\\
d=(-\sqrt{2}+1)e^{(2+\sqrt{2})t}+(\sqrt{2}+1)e^{(2-\sqrt{2})t}.
\end{array}
$$

\item Using formula (\ref{6.2.6}), we have
\begin{align}\label{6.2.12}
\begin{aligned}
x(t)&=x_0\lvplus x_1t\lvplus\int_0^{t}d\tau_1\int_0^{\tau_1}e^{A_*\tau_2}d\tau_2 x_2\\
    &=x_0\lvplus x_1t\lvplus\int_0^t(A_*^{-1}e^{A_*\tau_1}-A_*^{-1})d\tau_1 x_2\\
    &=x_0\lvplus x_1t\lvplus\left[A_*^{-2}e^{A_*t}-A_*^{-2}-A_*^{-1}t\right]x_2.\\
\end{aligned}
\end{align}
(\ref{6.2.12}) is exactly the same as (\ref{6.2.11}).
\end{itemize}
\end{enumerate}
\end{exa}

\section{Linear Control Systems}

\subsection{Discrete Time Time-Invariant Linear Control System}

Consider the following discrete time control system:
\begin{align}\label{7.1.1}
\begin{cases}
x(t+1)=A\lvtimes x(t)\lvplus B\ltimes u(t)\\
y(t)=C\lvtimes x(t),
\end{cases}
\end{align}
where $A\in {\cal M}_{m\times n}$,  $B\in {\cal M}_{m\times p}$, $C\in {\cal M}_{q\times m}$.

First, we consider the controllability.

\begin{dfn} \label{d7.1.1} System (\ref{7.1.1}) is
\begin{enumerate}
\item controllable from $x_0$ to $x_d$, where $x_0,~x_d \in {\cal V}$, if there exist a $T>0$ and a sequence $u_0,u_1,\cdots,u_{T-1}$, such that the trajectory can be driven from $x(0)=x_0$ to $x(T)=x_d$;
\item controllable from ${\cal V}_i$ to ${\cal V}_j$, if it is controllable from any $x_0\in {\cal V}_i$ to any $x_d\in {\cal V}_j$.
\end{enumerate}
\end{dfn}

For statement ease, the first (second) type of controllability is called the point (space) controllability. We are particularly interested in space controllability, because later on you will see that the controllability is independent on the particular starting and ending points. This is exactly the same as the standard linear system.

Assume $A$ is dimension-bounded. Then from any $x_0\in {\cal V}$, say, $x_0\in {\cal V}_{r_0}$, there exists a unique $A$-invariant space ${\cal V}_{r_*}$, depending on $r_0$, and assume $r_*=r_s$. Now restrict system (\ref{7.1.1}) on ${\cal V}_{r_*}$ we have
\begin{align}\label{7.1.2}
\begin{array}{l}
\begin{cases}
x(t+1)=A_sx(t)+B_su(t)\\
y(t)=C_sx(t),
\end{cases}\\
\mbox{where}\\
A_s=A_*,\; B_s=(B\otimes \J_{r_*/m}),\; C_s=\left(C\otimes I_{r_*/m}\right).
\end{array}
\end{align}
(\ref{7.1.2}) is a standard linear control system, which is called the stationary  realization of (\ref{7.1.1}). Note that the stationary realization depends on $x_0\in {\cal V}_{r_0}$.

Back to controllability. The following result is obvious.

\begin{thm}\label{t7.1.3}
\begin{enumerate}
\item System (\ref{7.1.1}) is controllable from ${\cal V}_i$ to ${\cal V}_j$, if
\begin{itemize}
\item[(i)] $A$ is dimension-bounded.
\item[(ii)] $r_0:=\dim({\cal V}_i)$ and $r_*:=\dim({\cal V}_j)$ satisfy (\ref{5.2.2}). Precisely speaking, starting from $r_0$, the algorithm (\ref{5.2.2}) terminates at $r_*$.
\item[(iii)] The stationary realization (\ref{7.1.2}) is controllable.
\end{itemize}
\item If {i}-(ii) hold, then (iii) is also necessary.
\end{enumerate}
\end{thm}

\noindent{\it Proof}.
\begin{enumerate}
\item First, we show the restriction is well defined. Note that since $A\in {\cal M}_{m\times n}$, and $A$ is dimension-bounded, then $n=m\mu_x$. It follows that $m|r_*$.
Then the restriction of (\ref{7.1.1}) on ${\cal V}_{r_*}$ becomes (\ref{7.1.2}),
%\begin{align}\label{7.1.101}
%\begin{array}{ccl}
%x(t+1)&=&A_*x(t)\lvplus Bu=A_*x(t)+(Bu)\otimes \J_{r_*/m}\\
%      &=&A_*x(t)+(B\otimes \J_{r_*/m})u,\\
%\end{array}
%\end{align}
which is a classical control system. Hence the controllability of (\ref{7.1.1}) is independent of the entrance position of the trajectory into ${\cal V}_{r_*}$.
Then the space controllability is obvious.
\item
Since the stationary realization (\ref{7.1.2}) is a standard control system and the trajectory of (\ref{7.1.1}) coincides with the one of (\ref{7.1.2}) after it reaches ${\cal V}_{r_*}$. The conclusion is obvious.
\end{enumerate}
\hfill $\Box$

Next, as we know that the trajectory starting from ${\cal V}_{r_0}$ will also reach ${\cal V}_{\ell}$, where $\ell\in\{r_0,r_1,\cdots,r_{s-1}\}$. ($r_s=r_*$ and ${\cal V}_{r_*}$ is $A$-invariant). Then is system (\ref{7.1.1}) controllable from ${\cal V}_{r_0}$ to ${\cal V}_{r_\ell}$?

Let $\ell|s$, say $s=\ell j$. Then \cite{chepr} we define an embedding mapping, $\bed :{\cal V}_{\ell}\ra {\cal V}_{s}$, as
$$
\bed_j(x)=x\otimes \J_j.
$$

\begin{prp}\label{7.1.4} Consider ${\cal V}_{r_k}$, where $0<k<s$ and $r_s= r_*$. Then the controllable subspace $C_{r_k}$ is
\begin{align}\label{7.1.3}
C_{r_k}=\bed_j\left(\Span\{B,A\lvtimes B,\cdots,A^{k-1}\lvtimes B\}\right).
\end{align}
where $j=\frac{r_k}{m}$.
\end{prp}
\noindent{\it Proof}. First, we claim that the trajectory reaches  ${\cal V}_{\ell}$ exactly once, where $\ell=r_k$. Assume the trajectory reaches ${\cal V}_{\ell}$ more than once. Then there is a cycle of spaces
$$
{\cal V}_{\ell}\ra {\cal V}_{\ell+1}\ra \cdots \ra {\cal V}_{\ell+T}= {\cal V}_{\ell}.
$$
Then the trajectory will never reach ${\cal V}_{r_*}$, which is a contradiction.

Calculating the reachable set yields that
$$
R_k=\left\{A^kx_0\lvplus Bu_{k-1}\lvplus A\lvtimes Bu_{k-2}\lvplus \cdots \lvplus A^{k-1}\lvtimes Bu_0\right\}
$$
Note that $\Span\{A\}\lvplus \Span\{B\}\subset {\cal V}_m$, the conclusion follows.
\hfill $\Box$

\begin{exa}\label{e7.1.4} Consider system (\ref{7.1.1}), where
$$
A=\begin{bmatrix}
1&2&-1&0\\
1&-2&2&-1
\end{bmatrix},
$$, and $B=[1,0]^T$.
\begin{enumerate}
\item Is it controllable from $\R^3$ to $\R^6$?

Then stationary realization is
\begin{align}\label{7.1.4}
x(t+1)=A_*x(t)+B_su,
\end{align}
where $A_*$ is obtained in (\ref{6.2.8}), and
$$
B_s=\begin{bmatrix}1\\0\end{bmatrix}\otimes \J_3.
$$
Then it is easy to verify that (\ref{7.1.4}) is not controllable. According to Theorem \ref{t7.1.3}, system (\ref{6.1.3}) is not controllable from $\R^3$ to $\R^6$.

\item Is it controllable from $\R^8$ to $\R^2$?

We also have the corresponding stationary realization (\ref{6.1.3}) with
$A_*$ as in (\ref{6.2.10}), and $B_s=B$. Then it is easy to verify that this stationary realization is controllable. So system (\ref{6.1.3}) is  controllable from $\R^8$ to $\R^2$.

\item The reachable subspace from $\R^8$ to $\R^4$?

Using (\ref{7.1.3}), the controllable subspace is:
$$
\begin{array}{ccl}
C_{r_1}&=&\bed_2\left(\Span\{B\}\right)\\
~&=&\Span\{[1,1,0,0]^T\}.
\end{array}
$$
\end{enumerate}
\end{exa}

Next we consider the observability of system (\ref{7.1.1}).  We assume $A$ is dimension-bounded.

\begin{dfn}\label{d7.1.5} Consider system (\ref{7.1.1}). The system is observable from $\R^{r_0}$, if for any  $x(0)\in \R^{r_0}$ there exists a $T>0$ such that $x(T)$ can be determined by the outputs $y(t)$, $t\geq 0$.
\end{dfn}

The following result is obvious.

\begin{thm}\label{t7.1.6} Consider system (\ref{7.1.1}) and assume $A$ is dimension-bounded. (\ref{7.1.1}) is observable, if and only if, its corresponding stationary  realization (\ref{7.1.2})
is observable.
\end{thm}

\begin{exa}\label{e7.1.7} Consider system (\ref{7.1.1}) with $A$ and $B=b$ as in Example \ref{e7.1.4}, and  $C=[0,1]$.
\begin{enumerate}
 \item Assume $x_0\in \R^3$. Then the $A$-invariant space is $\R^6$. The stationary realization on $\R^6$ has
 $A_s$ and $B_s$ as in 1) of Example \ref{e7.1.4}, and
 $$
 C_s=C\otimes I_3.
 $$
 Then it is ready to prove that its stationary realization is observable.
 \item Assume $x_0\in \R^8$. Then the $A$-invariant space is $\R^2$. The stationary realization on $\R^2$ has $A_s$ and $B_s$ as in 2) of Example \ref{e7.1.4}, and $C_s=C$. It is also observable.
\end{enumerate}
\end{exa}

\subsection{Continuous-time Case}

This section considers a continuous time control system:
\begin{align}\label{7.2.1}
\begin{cases}
\dot{x}(t)=A\lvtimes x(t)\lvplus Bu(t)\\
y(t)=C \lvtimes x(t),
\end{cases}
\end{align}
where $A\in {\cal M}_{m\times n}$,  $B\in {\cal M}_{m\times p}$, $C\in {\cal M}_{q\times m}$.

First, we calculate a trajectory of (\ref{7.2.1})

A straightforward computation shows the following:
\begin{align}\label{7.2.2}
A\lvtimes (Bu)=(A\lvtimes B)u=\dsum_{j=1}^pu_j\left(A\lvtimes \Col_j(B)\right).
\end{align}
Then we can express the trajectory of (\ref{7.2.1}) as follows:
\begin{align}\label{7.2.3}
x(t)=e^{At}\lvtimes x_0\lvplus\int_0^t e^{A(t-\tau)}\lvtimes Bu(\tau) d\tau.
\end{align}
(\ref{7.2.3}) can be proved by a straightforward computation.

We already know how to calculate the drift term of (\ref{7.2.3}). Now we calculate the integral part:
Similar to (\ref{6.2.3}) we can calculate that

\begin{align}\label{7.2.4}
\begin{array}{l}
e^{A(t-\tau)}\lvtimes Bu(\tau)\\
~=Bu_0\lvplus (t-\tau)A\lvtimes Bu(\tau)\lvplus\frac{(t-\tau)^2}{2!}A^2\lvtimes Bu(\tau)\lvplus\cdots\\
~~~\lvplus\frac{(t-\tau)^s}{s!}A^s\lvtimes Bu(\tau)\lvplus\cdots\\
~=\,\left[Bu(\tau)\lvplus (t-\tau)A\lvtimes Bu(\tau)\lvplus\frac{(t-\tau)^2}{2!}A^2\lvtimes Bu(\tau)\lvplus\cdots\right.\\
~~~\left.\lvplus\frac{(t-\tau)^{s-1}}{(s-1)!}A^{s-1}\lvtimes Bu(\tau)\right]\\
~~~\lvplus\left[\frac{(t-\tau)^s}{s!}I_{r_*}+\frac{t^{s+1}}{(s+1)!}A_*+\cdots \right]Bu(\tau)\\
~:=I\lvplus II,
\end{array}
\end{align}
where
\begin{align}\label{7.2.5}
\begin{array}{ccl}
I&=&\left[Bu(\tau)\lvplus (t-\tau)A\lvtimes Bu(\tau)\lvplus\frac{(t-\tau)^2}{2!}A^2\lvtimes Bu(\tau)\right.\\
~&~& \left.\lvplus\cdots\lvplus\frac{(t-\tau)^{s-1}}{(s-1)!}A^{s-1}\lvtimes Bu(\tau)\right]\\
II&=&\left[\frac{(t-\tau)^s}{s!}I_{r_*}+\frac{t^{s+1}}{(s+1)!}A_*+\cdots \right]B_s,
\end{array}
\end{align}
with
$$
B_s=A^{s}\lvtimes Bu(\tau).
$$
If $A_*$ is invertible, then similar to (\ref{6.2.4}) we have
\begin{align}\label{7.2.6}
\begin{array}{ccl}
II&=&A_*^{-s}\left( e^{ A_*  t}\lvminus \left[I_{r_*}+(t-\tau)A_*+\frac{t^2}{2!}A_*^2+\cdots\right.\right.\\
~&~&\left.\left.+\frac{(t-\tau)^{s-1}}{(s-1)!}A_*^{s-1}
\right]\right)B_s.
\end{array}
\end{align}
Otherwise, similar to (\ref{6.2.6}) we have
\begin{align}\label{7.2.7}
II=\int_0^td\tau_1\int_0^{\tau_1}d\tau_2\cdots\int_0^{\tau_{s-1}}e^{A_*(\tau_s-\tau)}d\tau_s B_s.
\end{align}
Finally,
\begin{align}\label{7.2.8}
\int_0^t e^{A(t-\tau)}\lvtimes Bu(\tau) d\tau=\int_0^t (I+II)  d\tau,
\end{align}
which is a classical integration.
We give a numerical description for this.

\begin{exa}\label{e7.2.1} Consider the following system:
\begin{align}
\label{7.2.9}
\dot{x}(t)=\begin{bmatrix}
1&0&-1&0\\
0&-1&0&1
\end{bmatrix}\lvtimes
x(t)\lvplus \begin{bmatrix}
1\\0\end{bmatrix}
u.
\end{align}
with $x(0)=x(0)=[1,0,1,-1]^T$.
It is easy to calculate that
$r_*=r_1=2$, $A\lvtimes x(0)=[0,-1]^T$, and
$$
A_0=A_*=\begin{bmatrix}
1&-1\\-1&1
\end{bmatrix}.
$$
Calculate the Jordan canonical form of $A_*$ as
$$
\tilde{A_*}:=P^{-1}A_*P=\begin{bmatrix}0&0\\0&2\end{bmatrix},
$$
where
$$
P=\begin{bmatrix}1&1\\1&-1\end{bmatrix}.
$$
Using formula (\ref{6.2.6}), we have
$$
\begin{aligned}
e^{At}x_0&=x_0\lvplus \int_0^t e^{A_*\tau}d\tau (A\lvtimes x_0)\\
 &=x_0\lvplus P\int_0^t e^{\tilde{A}_*}P^{-1}\begin{bmatrix}0\\-1\end{bmatrix}\\
&=x_0\lvplus \frac{1}{2}\begin{bmatrix}
t+\frac{1}{2}(e^{2t}-1)&t-\frac{1}{2}(e^{2t}-1)\\
t-\frac{1}{2}(e^{2t}-1)&t+\frac{1}{2}(e^{2t}-1)\\
\end{bmatrix} \begin{bmatrix}0\\-1\end{bmatrix}\\
 &=x_0\lvminus \frac{1}{2}\begin{bmatrix}
t-\frac{1}{2}(e^{2t}-1)\\
t+\frac{1}{2}(e^{2t}-1)\\
\end{bmatrix}.
\end{aligned}
$$
Since ${\cal V}_{r_*}=\R^2$, we have
$$
\begin{array}{l}
 e^{A(t-\tau)}\lvtimes Bu(\tau)\\
 =Bu(\tau)\lvplus \\
\frac{1}{2}\begin{bmatrix}
(t-\tau)+\frac{1}{2}(e^{2(t-\tau)}-1)&(t-\tau)-\frac{1}{2}(e^{2(t-\tau)}-1)\\
(t-\tau)-\frac{1}{2}(e^{2(t-\tau)}-1)&(t-\tau)+\frac{1}{2}(e^{2(t-\tau)}-1)\\
\end{bmatrix}\\
~\ltimes A\lvtimes Bu(\tau)\\
=\left[B+\frac{1}{2} \begin{bmatrix}
e^{2(t-\tau)}-1\\
-e^{2(t-\tau)}-1\\
\end{bmatrix}\right]u(\tau).
\end{array}
$$

We conclude that the trajectory of system (\ref{7.2.9}), starting from $x_0$, is
\begin{align}\label{7.2.10}
\begin{aligned}
x(t)=&\;
x_0\lvminus \frac{1}{2}\begin{bmatrix}
t-\frac{1}{2}(e^{2t}-1)\\
t+\frac{1}{2}(e^{2t}-1)\\
\end{bmatrix}\\
&\;\lvplus \int_0^t
\left[B+\frac{1}{2} \begin{bmatrix}
e^{2(t-\tau)}-1\\
-e^{2(t-\tau)}-1\\
\end{bmatrix}\right]u(\tau)d\tau.
\end{aligned}
\end{align}
\end{exa}

For continuous-time time-invariant linear system, the controllability and observability cannot be defined as for discrete time case.
We turn to its stationary realization.

\begin{prp}\label{p7.2.2} Consider system (\ref{7.2.1}) and assume $A$ is dimension-bounded. Then for any $x_0\in {\cal V}$ we have
a corresponding stationary  realization as
\begin{align}\label{7.2.11}
\begin{array}{l}
\begin{cases}
\dot{x}(t)=A_sx(t)+B_su,\\
y(t)=C_sx(t).
\end{cases}\\
\mbox{where}\\
A_s=A_*,\; B_s=\left(B\otimes \J_{r_*/m}\right),\;C_s=\left(C\otimes I_{r_*/m}\right).
\end{array}
\end{align}
\end{prp}

\noindent{\it Proof}. Since $A$ is dimension-bounded the existence of ${\cal V}_{r_*}$ is assured. From the calculation of trajectory one sees easily that the system will eventually evolve on ${\cal V}_{r_*}$, which depends on $x_0$. Then on this ${\cal V}_{r_*}$ the stationary realization becomes (\ref{7.2.11}).
\hfill $\Box$

Taking Theorems \ref{t7.1.3} and \ref{t7.1.6} into consideration, the controllability and observability of system (\ref{7.2.1}) can reasonably be defined as follows:

\begin{dfn}\label{d7.2.3}
\begin{enumerate}
\item The continuous time linear system (\ref{7.2.1}) is controllable, if its stationary realization (\ref{7.2.11}) is controllable.
\item The continuous time linear system (\ref{7.2.1}) is observable, if its stationary realization (\ref{7.2.11}) is observable.
\end{enumerate}
\end{dfn}

\begin{exa}\label{e7.2.4}
Consider the following system
\begin{align}\label{7.2.12}
\begin{cases}
\dot{x}(t)=[1,2,-1]\lvtimes x(t)+2u(t),\\
y(t)=-x(t).
\end{cases}
\end{align}
with $x_0=[1,-1]^T$. Check if the system with $x(0)=x_0$ is controllable? observable?

It is easy to calculate that $r_*=r_0=2$, and then
$$
A_{*}=\begin{bmatrix}
3&-1\\1&1\end{bmatrix},
$$
$$
B_s=B\otimes \J_2=[2,2]^T,
$$
and
$$
C_s=C\otimes I_2=-I_2.
$$
The stationary realization of (\ref{7.2.12}) is
\begin{align}\label{7.2.13}
\begin{cases}
\dot{x}(t)=\begin{bmatrix}3&-1\\1&1\end{bmatrix}x(t)+\begin{bmatrix}2\\2\end{bmatrix}u(t),\\
y(t)=\begin{bmatrix}-1&0\\0&-1\end{bmatrix}x(t).
\end{cases}
\end{align}
It is not controllable. It is observable.
\end{exa}

\section{Quotient Space}

\subsection{Quotient Vector Space $\Omega$}

Recall {\bf Case 2}  in Subsection III-A, where $x=-(-x)$ and assume
$z=-x$ is defined as $x+z\in \vec{0}$. Then it is easy to see that $y=-(z)$, if and only if, there exist $\J_{\a}$ and $\J_{\b}$ such that
\begin{align}\label{8.1.1}
x\otimes \J_{\a}=y \otimes \J_{\b}.
\end{align}

Similarly, if we check the distance between $x$ and $y$ using (\ref{3.2.5}), one sees easily that  $d(x,y)=0$, if and only if, (\ref{8.1.1}) holds.
These facts lead the following definition.

\begin{dfn}\label{d8.1.1} Let $x,~y\in {\cal V}$. Then $x$ and $y$ are said to be vector equivalent (V-equivalent), denoted by $x\lra y$, if there exist  $\J_{\a}$ and $\J_{\b}$ such that (\ref{8.1.1}) holds.
\end{dfn}

\begin{rem}\label{r8.1.2}
\begin{enumerate}
\item It is necessary to verify that $\lra$ is an equivalence relation. That is, (i) $x\lra x$; (ii) $x\lra y$, if and only if, $y\lra x$; (iii) $x\lra y$ and $y\lra z$ implies $x\lra z$. The verification is trivial.
\item The equivalence class of $x$ is denoted as
$$
\bar{x}:=\{y\;|\; y\lra x\}.
$$
\item The quotient space is defined as
$$
\Omega:={\cal V}/\lra=\{\bar{x}\;|\;x\in {\cal V}\}.
$$
\end{enumerate}
\end{rem}

We extend the V-addition to $\Omega$ as
\begin{align}\label{8.1.2}
\bar{x}\lvplus \bar{y}:=\overline{x\lvplus y}, \quad \bar{x},~\bar{y}\in \Omega.
\end{align}
Correspondingly,
\begin{align}\label{8.1.3}
\bar{x}\lvminus \bar{y}:=\bar{x}\lvplus (-\bar{y}), \quad \bar{x},~\bar{y}\in \Omega.
\end{align}

\begin{prp}\label{p8.1.3} \cite{chepr} (\ref{8.1.2}) and (\ref{8.1.3}) are properly defined. That is, if $x\lra x'$ and $y\lra y'$, that $x\lvplus y\lra x'\lvplus y'$.
\end{prp}

Define a scalar product $\cdot:\R\times \Omega\ra \Omega$ as
\begin{align}\label{8.1.4}
a\cdot \bar{x}:=\overline{ax}, \quad a\in \R,\; \bar{x}\in \Omega.
\end{align}

\begin{thm}\label{t8.1.4} $\Omega$ with addition $\lvplus$ and scalar product $\cdot$ defined by (\ref{8.1.4}) is a vector space.
\end{thm}

\noindent{\it Proof}. It follows from definition that the linearity of $\Omega$ comes from that of ${\cal V}$ directly. Moreover, since
$\vec{0}=\bar{0}\in \Omega$ is unique now, $-\bar{x}=\overline{-x}$ is also unique.
\hfill $\Box$

\subsection{Quotient Matrix Space $\Sigma$}

\begin{dfn}\label{d8.2.1} Let $A,~B\in {\cal M}$. $A$ and $B$ are said to be matrix equivalent (M-equivalent), denoted by $A\sim B$, if there exist $I_{\a}$ and $I_{\b}$ such that
\begin{align}\label{8.2.1}
A\otimes I_{\a}=B \otimes I_{\b}.
\end{align}
\end{dfn}

 The following facts have been proved in \cite{chepr}:

\begin{prp}\label{p8.2.2}
\begin{enumerate}
\item $\sim$ is an equivalence relation.
\item The equivalence class of $A$ is denoted as
$$
\A:=\{B\;|\;B\sim A\}.
$$
\item The quotient space is defined as
$$
\Sigma:={\cal M}/\sim.
$$
\item The operator $\ltimes: \Sigma\times \Sigma\ra \Sigma$ is defined as
\begin{align}\label{8.2.2}
\A\ltimes \B:=\left<A\ltimes B\right>.
\end{align}
Moreover, (\ref{8.2.2}) is properly defined.
\end{enumerate}
\end{prp}

As an immediate consequence of (\ref{8.2.2}), we have
\begin{cor}\label{c8.2.3} $(\Sigma,~\ltimes)$ is a monoid.
\end{cor}

Next, we consider the subset ${\cal M}_{\mu}$, where $\mu\in \Q_+$. Define
$$
\Sigma_{\mu}:={\cal M}_{\mu}/\sim.
$$
Then similar to vector case we have
\begin{prp}\label{p8.2.4} \cite{chepr}
\begin{enumerate}
\item Let $A\in {\cal M}_{m\times n}\subset {\cal M}_{\mu}$ and $B\in {\cal M}_{p\times q}\subset {\cal M}_{\mu}$. Then
\begin{align}\label{8.2.3}
\A\lplus \B:=\left<A \lplus B\right>.
\end{align}
(\ref{8.2.3}) is properly defined.

Correspondingly,
\begin{align}\label{8.2.301}
\A\lminus \B:=\left<A \lminus B\right>.
\end{align}

\item The scalar product is defined as
\begin{align}\label{8.2.4}
a\A:=\left<aA\right>,\quad a\in \R,\; A\in {\cal M}.
\end{align}
\item $\Sigma_{\mu}$ with M-addition $\lplus$ defined by (\ref{8.2.3}) and scalar product defined by (\ref{8.2.4}), $\left(\Sigma_{\mu},\lplus,\cdot\right)$ is a vector space.
\end{enumerate}
\end{prp}

\subsection{Dynamic System $(\Sigma,\lvtimes,\Omega)$}

The action of $\Sigma$ on $\Omega$ is defined as
\begin{align}\label{8.3.1}
\A\lvtimes \bar{x}:=\overline{A\lvtimes x}.
\end{align}
Then one can show easily that (\ref{8.3.1}) is properly defined. Moreover, we have

\begin{lem}\label{l8.3.1} $(\Sigma,\lvtimes,\Omega)$ is an S-system.
\end{lem}

\noindent{\it Proof}. The only thing needs to be proved is:
\begin{align}\label{8.3.2}
\A\lvtimes (\B\lvtimes \bar{x})=(\A \ltimes \B)\lvtimes \bar{x}.
\end{align}
Note that
$$
\begin{array}{l}
\A\lvtimes (\B\lvtimes \bar{x})=\A \lvtimes (\overline{B\lvtimes x})\\
~=\overline{A\lvtimes (B\lvtimes x)}=\overline{(A\ltimes B)\lvtimes x}\\
~=(\left<A\ltimes B\right>)\lvtimes \bar{x}=(\A\ltimes \B)\lvtimes \bar{x}.
\end{array}
$$
(\ref{8.3.1}) is proved.
\hfill $\Box$

 Define
\begin{align}\label{8.3.3}
\|\bar{x}\|_{{\cal V}}:=\|x\|_{{\cal V}},\quad \bar{x}\in \Omega.
\end{align}
\begin{align}\label{8.3.4}
\|\A\|_{{\cal V}}:=\|A\|_{{\cal V}},\quad \A\in \Sigma.
\end{align}
Then
\begin{lem}\label{l8.3.2} Both (\ref{8.3.3}) and (\ref{8.3.4}) are well defined.
\end{lem}

\noindent{\it Proof}.
\begin{itemize}
\item Consider (\ref{8.3.3}). It is enough to show that if $x\lra y$ then $\|x\|_{{\cal V}}=\|y\|_{{\cal V}}$. Assume $x\in {\cal V}_m$ and $y\in {\cal V}_n$. Then (refer to \cite{chepr}) there exists a $z\in {\cal V}$, called the largest common divisor of $x$ and $y$, such that
    \begin{align}\label{8.3.5}
    x=z\otimes \J_{\a},\quad y=z\otimes \J_{\b}.
    \end{align}
    It follows from (\ref{8.3.5}) that
    $$
    \frac{m}{\a}=\frac{n}{\b}:=r.
$$
Using (\ref{3.2.2}), we have
$$
\begin{array}{ccl}
\|x\|_{{\cal V}}&=&\sqrt{\frac{1}{m}\left<z\otimes \J_{\a}, z\otimes \J_{\a}\right>}\\
~&=&\sqrt{\frac{\a}{m}\left<z, z\right>}\\
~&=&\sqrt{\frac{1}{r}\left<z, z\right>}\\
\end{array}
$$
Similarly, we have
$$
\|y\|_{{\cal V}}=\sqrt{\frac{1}{r}\left<z, z\right>}.
$$
Hence (\ref{8.3.3}) is well defined.

\item  Consider (\ref{8.3.4}). It is enough to show that if $A\sim B$ then $\|A\|_{{\cal V}}=\|B\|_{{\cal V}}$. Assume $A\in {\cal M}_{m\times n}$ and $B\in {\cal M}_{p\times q}$. Then (refer to \cite{chepr}) there exists a $\Lambda\in {\cal M}$, called the largest divisor of $A$ and $B$, such that
    \begin{align}\label{8.3.6}
    A=\Lambda \otimes I_{\a},\quad B=\Lambda\otimes I_{\b}.
    \end{align}
    It follows from (\ref{8.3.6}) that
    $$
    \frac{m}{\a}=\frac{p}{\b},\quad
    \frac{n}{\a}=\frac{q}{\b}.
$$
Then we have
$$
 \frac{n}{m}=\frac{q}{p}:=r\\
$$
Using (\ref{3.2.4}), we have
$$
\begin{array}{ccl}
\|A \|_{{\cal V}}&=&\sqrt{\frac{n}{m}} \sqrt{\sigma_{\max}[(\Lambda\otimes I_{\a})^T (\Lambda\otimes I_{\a})]}\\
~&=&\sqrt{\frac{n}{m}} \sqrt{\sigma_{\max}[(\Lambda^T\Lambda)\otimes I_{\a})]}\\
~&=&\sqrt{r\sigma_{\max}(\Lambda^T\Lambda)}.
\end{array}
$$
Similarly, we have
$$
\begin{array}{ccl}
\|B\|_{{\cal V}}&=&\sqrt{\frac{q}{p}} \sqrt{\sigma_{\max}(\Lambda^T\Lambda)}\\
~&=&\sqrt{r\sigma_{\max}(\Lambda^T\Lambda)}.
\end{array}
$$
Hence (\ref{8.3.4}) is well defined.
\end{itemize}
\hfill $\Box$

\begin{thm}\label{t8.3.3} $(\Sigma,\lvtimes,\Omega)$ is a linear dynamic system.
\end{thm}

\noindent{\it Proof}. First, using the norm defined (\ref{8.3.3}) we have
\begin{align}\label{8.3.7}
d(\bar{x},\bar{y})=\|\bar{x}\lvminus \bar{y}\|_{{\cal V}}
=\|x\lvminus y\|_{{\cal V}}^2=d(x,y).
\end{align}
Since $\Omega$ is a vector space, using this distance $\Omega$ becomes a metric space. Since a metric space is $T_4$, which implies $T_2$ (i.e., Hausdorff)~\cite{kel75}, we know that $\Omega$ is Hausdorff.

Since
$$
\begin{array}{ccl}
\|\A\lvtimes \bar{x}\|_{{\cal V}}&=&\|A\lvtimes x\|_{{\cal V}}\\
~&\leq& \|A\|_{{\cal V}}\|x\|_{{\cal V}}\\
~&=&\|\A\|_{{\cal V}}\|\bar{x}\|_{{\cal V}},\quad \bar{x}\in \Omega.
\end{array}
$$
The mapping $\lvtimes:\Sigma\times \Omega\ra \Omega$ is continuous. Since the linearity is obvious, the conclusion follows.
\hfill $\Box$

\subsection{Formal Polynomial ${\cal P}_{\Sigma}$}

Consider the set of formal polynomials on quotient space $\Sigma$, which is defined as
\begin{align}\label{8.4.1}
{\cal P}_{\Sigma}:=\oplus_{\mu\in \Q_+}\Sigma_{\mu}.
\end{align}
As aforementioned that $\Sigma_{\mu}$, $\mu\in \Q_+$, are vector space. According to Proposition \ref{p4.1.2}, ${\cal P}_{\Sigma}$ is a vector space.
Let $\bar{p}\in {\cal P}_{\Sigma}$. Then $\bar{p}$ can be expressed into a formal polynomial as
\begin{align}\label{8.4.2}
\bar{p}=\dsum_{\mu\in \Q_+}\A_{\mu} z_{\mu},
\end{align}
where $\A_{\mu}\in \Sigma_{\mu}$.
$p\in {\cal P}_{{\cal M}}$ is called a representative of $\bar{p}$, if
\begin{align}\label{8.4.3}
p=\dsum_{\mu\in \Q_+}A_{\mu} z_{\mu},\quad A_{\mu}\in \A_{\mu}
\end{align}
As for ${\cal P}_{{\cal M}}$, we also assume {\bf A-1} for ${\cal P}_{\Sigma}$. That is

\vskip 2mm

\noindent{\bf A-2}
$$
\dsum_{\mu\in Q_+}\|\A_{\mu}\|_{{\cal V}}<\infty.
$$

\vskip 2mm

\begin{dfn}\label{d8.4.1}
\begin{enumerate}
\item Assume $\bar{p}\in {\cal P}_{\Sigma}$ is as in (\ref{8.4.2}). The action of ${\cal P}_{\Sigma}$ on $\Omega$ is defined as
\begin{align}\label{8.4.4}
\bar{p}\lvtimes \bar{x}:=\dsum_{\mu\in \Q_+} \A_{\mu}\lvtimes \bar{x},\quad \bar{x}\in \Omega.
\end{align}
\item The product on ${\cal P}_{\Sigma}$ is defined as
\begin{align}\label{8.4.5}
\bar{p}\ltimes \bar{q}:=\overline{p\ltimes q},\quad p\in \bar{p},\; q\in \bar{q}.
\end{align}
\end{enumerate}
\end{dfn}

The following result is straightforward verifiable:

\begin{prp}\label{p8.4.2}
\begin{enumerate}
\item Assume {\bf A-2}, then $\bar{p}$ is a proper operator on $\Omega$.
\item Assume both $\bar{p}$ and $\bar{q}$ satisfy {\bf A-2}, then so is $\bar{p}\ltimes \bar{q}$.
\end{enumerate}
\end{prp}

Hereafter, we assume {\bf A-2} is always true. In fact, we consider only $\bar{p}\in {\cal P}_{\Sigma}$, which satisfy {\bf A-2}. For notational ease, we still use ${\cal P}_{\Sigma}$ for this subset. Then using Theorem \ref{t8.3.3}, we have

\begin{cor}\label{c8.4.3}
$\left({\cal P}_{\Sigma},\lvtimes,\Omega\right)$ is a linear dynamic system.
\end{cor}

Corresponding to ${\cal P}_{{\cal M}^{\mu}}\subset {\cal P}_{{\cal M}}$,  ${\cal P}_{\Sigma^{\mu}}\subset {\cal P}_{\Sigma}$, is particularly important.
$$
{\cal P}_{\Sigma^{\mu}}=\oplus_{n=0}^{\infty}{\cal P}_{\Sigma_{{\mu}^n}}.
$$
$\bar{p}\in {\cal P}_{\Sigma^{\mu}}$ can be expressed in a formal polynomial as
\begin{align}\label{8.4.6}
\bar{p}=\dsum_{n=0}^{\infty}\A_n z^n,
\end{align}
where $\A_n\in \Sigma_{{\mu}^n}$.
Particularly, assume $f(x)$ is an analytic function with its Taylor expansion as in (\ref{4.2.8}). Then
$$
f(\A):=\dsum_{n=0}^{\infty} c_n\A^n z^n,
$$
and
$$
f(\A)\lvtimes \bar{x}=\dsum_{n=0}^{\infty} c_n\overline{A^n\lvtimes x}.
$$

Finally, we introduce a special kind of formal polynomials: If $\A_i=c_i\left<A^i\right>$, $i=0,1,\cdots$, for a fixed $A$, then it is described as
\begin{align}\label{8.4.7}
{\cal P}(A):=\left\{\dsum_{i=0}^{\infty} c_i \left<A^i\right> z^{i}\;|\; c_i\in \R\right\}.
\end{align}
$p\in {\cal P}(A)$ is called a principle formal polynomial (PFP). In the sequel, we will see that PFP is particular important in investigating time invariant linear systems.

\subsection{Lie Algebra of Formal Polynomial ${\cal P}_{\Sigma}$}

Observe ${\cal P}_{\Sigma}$ again. One sees easily that
\begin{itemize}
\item ${\cal P}_{\Sigma}$ with $\lplus$ is a vector space.
\item ${\cal P}_{\Sigma}$ with $\ltimes$ is a monoid.
\end{itemize}

Then can we pose a Lie algebraic structure on it? The answer is ``Yes". We refer to \cite{boo86,var84} for basic concepts of Lie algebra used in this subsection.

\begin{prp}\label{p8.5.1} Let $\bar{p},~\bar{q}\in {\cal P}_{\Sigma}$. Define a Lie bracket on ${\cal P}_{\Sigma}$ as
\begin{align}\label{8.5.1}
[\bar{p},~\bar{q}]:=\bar{p}\ltimes \bar{q} \lminus \bar{q}\ltimes \bar{p}.
\end{align}
Then ${\cal P}_{\Sigma}$ becomes a Lie algebra. Precisely speaking, the following are satisfied:
\begin{itemize}
\item[(1)] (Linearity)
\begin{align}\label{8.5.2}
[\bar{p}_1+\bar{p}_2,\bar{q}]=[\bar{p}_1,\bar{q}]\lplus [\bar{p}_2,\bar{q}];
\end{align}
\item[(2)] (Skew-Symmetry)
\begin{align}\label{8.5.3}
[\bar{p},\bar{q}]=-[\bar{q},\bar{p}];
\end{align}
\item[(3)] (Jacobi Identity)
\begin{align}\label{8.5.4}
[\bar{p},[\bar{q},\bar{r}]]\lplus [\bar{q},[\bar{r},\bar{p}]]\lplus [\bar{r},[\bar{p},\bar{q}]]=\bar{0},
\end{align}
where $\bar{0}$ is the zero element in ${\cal P}_{\Sigma}$.
\end{itemize}
\end{prp}

Let
$${\cal P}_{\Lambda}:=\{\bar{p}_{\lambda}\;|\;\lambda\in \Lambda\}\subset {\cal P}_{\Sigma}.
$$
Then we define a Lie algebra ${\cal P}_{\Lambda}^{LA}$, the smallest Lie subalgebra of ${\cal P}_{\Sigma}$ containing ${\cal P}_{\Lambda}$.

\begin{exa}\label{e8.5.2}
\begin{enumerate}
\item ${\cal P}_{\mu}$ is a Lie sub-algebra generated by $\Sigma_{\mu}$. That is,
$$
{\cal P}_{\mu}={\cal P}_{\Sigma_{\mu}}^{LA}.
$$
\item ${\cal P}(\A)$ is a Lie sub-algebra generated by $\A$. In fact, it is easy to prove that
$$
{\cal P}(\A)^{LA}={\cal P}(\A).
$$
That is, ${\cal P}(\A)$ itself is a Lie algebra.
\end{enumerate}
\end{exa}

\begin{rem}\label{r8.5.3} The most important subalgebra of ${\cal P}_{\Sigma}$ is ${\cal P}_1$, That is, $\{\A\;|\; A ~\mbox{is square}\}$. This subalgebra has very nice properties. We refer to \cite{chepr} for details.
\end{rem}

\section{Linear System on Quotient Space}

Assume we have linear systems on quotient space as follows:
\begin{itemize}
\item Discrete Time Linear System:
\begin{align}\label{9.1.1}
\begin{array}{l}
\bar{x}(t+1)=\A(t)\lvtimes \bar{x}(t),\quad \bar{x}(0)=\overline{x_0}\\
~~\bar{x}(t)\in \Omega,\;\A(t)\in \Sigma.
\end{array}
\end{align}
\item Continuous Time Linear System:
\begin{align}\label{9.1.2}
\begin{array}{l}
\dot{\bar{x}}(t)=\A(t)\lvtimes \bar{x}(t),\quad \bar{x}(0)=\overline{x_0}\\
~~\bar{x}(t)\in \Omega,\;\A(t)\in \Sigma.
\end{array}
\end{align}
\item Discrete Time Linear Control System:
\begin{align}\label{9.1.3}
\begin{cases}
\bar{x}(t+1)=\A(t)\lvtimes \bar{x}(t)\lvplus \bar{B}u(t),\quad \bar{x}(0)=\overline{x_0}\\
\bar{y}=\C \lvtimes\bar{x},
\end{cases}
\end{align}
where $A\in {\cal M}_{m\times n}$, $B\in {\cal M}_{q\times m}$, $C\in {\cal M}_{q\times m}$,  $\bar{x}(t)\in \Omega$, $\A(t)\in \Sigma$.
\item Continuous Time Linear Control System:
\begin{align}\label{9.1.4}
\begin{cases}
\dot{\bar{x}}(t)=\A(t)\lvtimes \bar{x}(t)\lvplus \bar{B}u(t),\quad \bar{x}(0)=\overline{x_0}\\
\bar{y}=\C \lvtimes\bar{x},
\end{cases}
\end{align}
where $A\in {\cal M}_{m\times n}$, $B\in {\cal M}_{q\times m}$, $C\in {\cal M}_{q\times m}$,  $\bar{x}(t)\in \Omega$, $\A(t)\in \Sigma$.
\end{itemize}

Then we construct the corresponding systems on original space as
\begin{itemize}
\item Discrete Time Linear System:
\begin{align}\label{9.1.5}
\begin{array}{l}
x(t+1)=A(t)\lvtimes x(t),\quad x(0)=x_0\\
~~x(t)\in {\cal V},\; A(t)\in {\cal M}.
\end{array}
\end{align}
\item Continuous Time Linear System:
\begin{align}\label{9.1.6}
\begin{array}{l}
\dot{x}(t)=A(t)\lvtimes x(t),\quad x(0)=x_0\\
~~x(t)\in {\cal V},\;A(t)\in {\cal M}.
\end{array}
\end{align}
\item Discrete Time Linear Control System:
\begin{align}\label{9.1.7}
\begin{cases}
x(t+1)=A(t)\lvtimes x(t)\lvplus Bu(t),\quad x(0)=x_0\\
y(t)=C\lvtimes x(t),
\end{cases}
\end{align}
where $A(t)\in {\cal M}_{m\times n}$, $B\in {\cal M}_{q\times m}$, $C\in {\cal M}_{q\times m}$,  $x(t)\in {\cal V}$.
\item Continuous Time Linear Control System:
\begin{align}\label{9.1.8}
\begin{cases}
\dot{x}(t)=A(t)\lvtimes x(t)\lvplus Bu(t),\quad x(0)=x_0\\
y=\C\lvtimes x,
\end{cases}
\end{align}
where $A\in {\cal M}_{m\times n}$, $B\in {\cal M}_{q\times m}$, $C\in {\cal M}_{q\times m}$,$x(t)\in {\cal V}$.
\end{itemize}

\begin{dfn}\label{d9.1.1} If there
are $A(t)\in \A(t)$, $B(t)\in \bar{B}(t)$, $C(t)\in \C(t)$, and $x_0\in \overline{x_0}$ such that
the systems (\ref{9.1.5})-(\ref{9.1.8}) exist, then   (\ref{9.1.5})-(\ref{9.1.8})  are called the realizations of  (\ref{9.1.1})-(\ref{9.1.4}) respectively.
\end{dfn}

The following result comes from (\ref{8.3.1}):

\begin{prp}\label{p9.1.2} Assume system (\ref{9.1.1}) (or (\ref{9.1.2}),  (\ref{9.1.3}), and (\ref{9.1.4})) has its realization (\ref{9.1.5}) (or (\ref{9.1.6}),  (\ref{9.1.7}), and (\ref{9.1.4}) respectively), and the trajectory of  (\ref{9.1.5}) (or (\ref{9.1.6}),  (\ref{9.1.7}),(\ref{9.1.8}))
is $x(t)$ with $x(0)=x_0$. Then $\bar{x}(t)=\overline{x(t)}$ is the trajectory of (\ref{9.1.1}) (or (\ref{9.1.2}),  (\ref{9.1.3}), and (\ref{9.1.4}) respectively) with $\bar{x}(0)=\overline{x_0}$.
\end{prp}

According to Proposition \ref{p9.1.2}, if a linear system defined on quotient space has its realization, then it is easy to find its trajectories. Hence the method proposed for linear pseudo dynamic systems based on $\left({\cal M},\lvtimes, {\cal V}\right)$  can be used directly for solving the linear dynamic systems based on $\left(\Sigma,\lvtimes, \Omega \right)$.

\begin{rem}\label{d9.1.3} The realization systems (\ref{9.1.5})-(\ref{9.1.8}) seem to be more ``realistic" than their parent systems  (\ref{9.1.1})-(\ref{9.1.4}), because the realization systems are evolving on real state space ${\cal V}$. But the realization systems have their own weakness, because their fundamental evolutionary model $({\cal M}, \lvtimes,{\cal V})$ is a pseudo dynamic system. While, the original systems on quotient spaces have their fundamental evolutionary model $(\Sigma,\lvtimes, \Omega)$, which is a dynamic system. Moreover, both $\Sigma$ and $\Omega$ are vector space.
\end{rem}

To further exploring the realization of systems on quotient space, we need some new concepts.

Define
$$
\overline{{\cal V}_r}:=\left\{\bar{x}\;|\;x\in {\cal V}_r\right\}.
$$

\begin{dfn}\label{d9.1.4} $\overline{{\cal V}_r}$ is called $\A$-invariant space, if
\begin{align}\label{9.1.9}
\A \lvtimes \overline{{\cal V}_r}\subset \overline{{\cal V}_r}.
\end{align}
\end{dfn}

Again  from (\ref{8.3.1}) we have the following result:

\begin{prp}\label{p9.1.5} If ${\cal V}_r$ is $A$-invariant, then $\overline{{\cal V}_r}$ is $\A$-invariant.
\end{prp}

A natural question is whether its converse is also true? Precisely speaking, if $\overline{{\cal V}_r}$ is $\A$-invariant, do there always exist a ${\cal V}_{r_s}\in \overline{{\cal V}_r}$ and an $A_s\in \A$, such that ${\cal V}_{r_s}$ is $A_s$-invariant? The answer is ``No". We give the following counter-example.
\begin{exa}\label{e9.1.6} Consider $\overline{{\cal V}_2}$ and $\A$, where
$$
A=\begin{bmatrix}a_{11}&0&a_{13}\\
a_{21}&0&a_{23}
\end{bmatrix}.
$$
Then for any $x=[\a,\b]^T\in {\cal V}_2$, we have
$$
\begin{array}{ccl}
\A\lvtimes\bar{x}&=&\overline{A\lvtimes x}\\
~&=&\overline{
\begin{bmatrix}
a_{11}&a_{12}\\
a_{21}&a_{22}\end{bmatrix}
\begin{bmatrix}
\a\\\b
\end{bmatrix} \otimes I_2}\in \overline{{\cal V}_2}.
\end{array}
$$
Hence $\overline{{\cal V}_2}$ is $\A$-invariant.

Now since $A\in {\cal M}_{2\times 3}$, any $\tilde{A}\in \A$ has $\mu_y=2$. According to Proposition \ref{p5.2.4}, $\tilde{A}$ has no invariant space.
\end{exa}

In fact, this fact makes the systems on quotient spaces more attractive.

\subsection{Projective Stationary Realization}

Example \ref{e9.1.6} shows that $\A$-invariant space $\bar{\cal V}_r$ may not comes from $A$-invariant space ${\cal V}_{s}$. As aforementioned that this phenomenon provides more geometric varieties for systems on quotient space. Even more, we may consider a subspace of ${\cal V}_r$. Let ${\cal S}_r\subset {\cal V}_r$ be a subspace. Define
$$
\bar{\cal S}_r:=\{y\;|\;\mbox{there exists}~ x\in {\cal S}_r,~mbox{such that}~y\lra x\}.
$$

\begin{dfn}\label{d9.2.1} Let ${\cal S}_r\subset {\cal V}_r$ be a subspace. $\bar{\cal S}_r$ is said to be $\A$-invariant, if
\begin{align}\label{9.2.1}
\A\lvtimes \bar{\cal S}_r\subset \bar{\cal S}_r.
\end{align}
\end{dfn}

As we discussed before, if there is an $x\neq 0$ in ${\cal V}_r$ and $A\lvtimes x\in {\cal V}_r$ then ${\cal V}_r$ is $A$-invariant. So if ${\cal S}_r\subset {\cal V}_r$ is $A$-invariant, then the whole space ${\cal V}_r$ is $A$-invariant. Hence, the $A$-invariant subspace ${\cal S}_r$ does not make sense. The following example shows that for quotient spaces  $\A$-invariant subspace $\bar{\cal S}_r$ does exit.

\begin{exa}\label{e9.2.2} Assume
$$
A=\begin{bmatrix}
a_{11}&a_{12}\\
a_{21}&a_{22}\\
a_{21}&a_{22}\\
\end{bmatrix},\quad
S_3=\Span\{\J_3\}\in {\cal V}_3.
$$
Then it is easy to verify that if
$$
a_{11}+a_{12}=a_{21}+a_{22}=a_{31}+a_{32},
$$
then
$\bar{S}_3$ is $\A$-invariant. It is well known that $A$ has no invariant space.
\end{exa}

 Consider a dynamic system
\begin{align}\label{9.2.2}
\bar{x}(t+1)=\A\lvtimes \bar{x}(t),\quad \bar{x}_0=\overline{x_0}.
\end{align}
Assume $\bar{x}_0\in \bar{\cal V}_r$ (or $\bar{x}_0\in \bar{\cal S}_r\subset \bar{\cal V}_r$), and
$\bar{\cal V}_r$ (correspondingly, $\bar{\cal S}_r$) is $\A$-invariant.
Then the trajectory $\bar{x}(t)\in \bar{{\cal V}}_r$ (or $\bar{x}(t)\in \bar{{\cal S}}_r$).
Then there exists a unique $x(t)\in {\cal V}_r$, such that
\begin{align}\label{9.2.3}
x(t)\in \bar{x}(t).
\end{align}
Define a projection
$$
Pr_r(\bar{x}(t)):=x(t),\quad x(t)\in {\cal V}_r\bigcup \bar{x}(t).
$$

\begin{dfn}\label{d9.2.3} The system
\begin{align}\label{9.2.4}
x(t+1)=A_*\lvtimes x(t),\quad x(0)=Pr_r(\bar{x}_0),
\end{align}
is called a projective stationary realization of (\ref{9.2.2}), if
\begin{align}\label{9.2.5}
x(t)=Pr_r(\bar{x}(t)),\quad t=1,2,\cdots.
\end{align}
\end{dfn}

Then we cane verify the following result directly.

\begin{prp}\label{p9.2.4}
Let $\bar{x}(t)$ be the trajectory of system (\ref{9.2.2}) with $\bar{x}_0=\overline{x_0}$. Moreover, $s_0$ is the smallest $s>0$ such that $\bar{x}(s)\in \bar{\cal V}_r$, and (\ref{9.2.4}) be a projective stationary realization of (\ref{9.2.2}) with $x_0=Pr_r(\bar{x}_{s_0})$.
Then after $s_0$ we have
\begin{align}\label{9.2.6}
Pr(\bar{x}(s_0+i))=x(i),\quad i=0,1,2,\cdots.
\end{align}
\end{prp}

A dynamic system over quotient spaces may have no stationary realization but with projective stationary realization. We give an example to show this.

\begin{exa}\label{e9.2.5}
\begin{enumerate}
\item Consider a discrete-time linear dynamic system as
\begin{align}\label{9.2.7}
\bar{x}(t+1)=\A\lvtimes \bar{x}(t),\quad \bar{x}_0=\overline{x_0},
\end{align}
where
$$
A=\begin{bmatrix}
1&2\\
1&2\\
-1&1\\
-1&1
\end{bmatrix};\quad x_0=\begin{bmatrix}
-1\\1\end{bmatrix}
$$
It is easy to see that $\bar{\cal V}_2$ is an invariant space. Because let $x=(\a,\b)^T\in {\cal V}_2$. Then
$$
A\lvtimes x=
\begin{bmatrix}
\a+2\b\\
\a+2\b\\
\b-\a\\
\b-\a
\end{bmatrix}=\begin{bmatrix}
\a+2\b\\
\b-\a\\
\end{bmatrix}\otimes \J_2.
$$
Then the system (\ref{9.2.7}) has a projective stationary realization as
$$
x(t+1)=A_*\lvtimes x(t)=
\begin{bmatrix}
1&2\\
-1&1\\
\end{bmatrix}\lvtimes x(t).
$$
Hence the trajectory of (\ref{9.2.7}) is
$$
\begin{array}{ccl}
\bar{x}(t)&=&\left< A_*\right>\lvtimes \bar{x}_0\\
~&=&\overline{
\begin{bmatrix}
1&2\\
-1&1\\
\end{bmatrix}^tx_0}.
\end{array}
$$
\item Consider a continuous-time linear dynamic system as
\begin{align}\label{9.2.8}
\dot{\bar{x}}(t)=\A\lvtimes \bar{x}(t),\quad \bar{x}_0=\overline{x_0},
\end{align}
where $A$ and $x_0$ are  as in 1).

Then it is easy to figure out that the trajectory is
$$
\bar{x}(t)=\overline{e^{A_*t}x_0}
$$
where as in 1) we have
$$
A_*=\begin{bmatrix}
1&2\\
-1&1\\
\end{bmatrix}
$$
\item Consider a continuous-time linear control system
\begin{align}\label{9.2.9}
\begin{cases}
\dot{\bar{x}}(t)=\A\lvtimes \bar{x}(t)\lvplus \bar{B}u, \\
\bar{y}(t)=\C\lvtimes \bar{x}(t),
\end{cases}
\end{align}
where $\bar{x}_0=\overline{x_0}$, $A$ is as before, and
$$
B=\begin{bmatrix}
1\\
2\\
\end{bmatrix},\quad
C=[1,0].
$$

It is easy to calculate that the trajectory is
$$
\bar{x}(t)=\overline{e^{A_*t}x_0+\int_{0}^te^{A_*(t-\tau)}Bu(\tau)d\tau}.
$$
Output is
$$
\bar{y}(t)=\overline{C\lvtimes x(t)}.
$$

Then it is ready to check that system (\ref{9.2.9}) is controllable and observable.

\end{enumerate}
\end{exa}

\section{Conclusion}

In this paper the cross-dimensional linear system is proposed and investigated. Roughly speaking, the linear system is constructed as follows: In the light of M-product (semi-tensor product) the set of matrices with arbitrary dimensions, ${\cal M}$, becomes a monoid. Consider the set of vectors with arbitrary dimensions, ${\cal V}$ as the objective state space. Then the V-product is defined as the action of ${\cal M}$ on ${\cal V}$, which provides a dynamic system, called the discrete-time S-system.

Two matrix products are considered. M-product is used for composing linear mappings; V-product is used for realizing linear mappings over object state space. Two matrix additions are proposed. M-addition makes set of matrices with different dimensions, ${\cal M}_{\mu}$, a hybrid vector space, and V-addition makes vectors of different dimensions, ${\cal V}$, a hybrid vector space.

Two matrix products and two matrix/vector additions are all the generalizations of conventional matrix product and matrix/vector addition. Inspired by them, the matrix equivalence and vector equivalence become natural.

Two equivalences  yield two quotient spaces. M-equivalence, $\sim$, is an equivalence relation on matrices ${\cal M}$, which makes the quotient space $\Sigma_{\mu}={\cal M}_{\mu}/\sim$ a vector space. V-equivalence, $\lra$, is an equivalence relation on vectors ${\cal V}$, which makes the quotient space $\Omega={\cal V}/\lra$ a vector space too.

The above concepts/properties form a set of new results of new matrix theory, where the dimension barrier of traditional matrix theory has been removed. It may be called the dimension-free matrix theory. As an application of these new concepts/operators the cross-dimensional linear dynamic systems over the quotient spaces are investigated in detail.

Using direct sum of vector spaces to combine all ${\cal M}_{\mu}$ ($\Sigma_{\mu}$) together yields the formal polynomial ${\cal P}_{{\cal M}}$ (correspondingly, ${\cal P}_{\Sigma})$, which is a hybrid vector space (correspondingly, vector space) over all matrices.

 Using formal polynomial, we consider the action of ${\cal M}_{{\cal M}}$ on ${\cal V}$ (${\cal P}_{\Sigma}$ on $\Omega$), which is considered  as the action of any matrices on any vectors, is a pseudo-dynamic system, $\left({\cal P}_{{\cal M}},\lvtimes, {\cal V}\right)$ (correspondingly, dynamic system $\left({\cal P}_{\Sigma},\lvtimes, \Omega\right)$),  where both ${\cal P}_{{\cal M}}$ and ${\cal V}$ are hybrid vector spaces (correspondingly, ${\cal P}_{\Sigma}$ and $\Omega$ are vector spaces), and the action is linear. This pseudo-dynamic system (correspondingly, dynamic system) is the cross-dimensional hybrid linear system (correspondingly, linear system).

Both discrete time and continuous time cross-dimensional (pseudo-) linear systems are investigated. Some properties are obtained. Particularly, as $\A$ is dimension-bounded the trajectory of the system is finitely computable.
In addition to dynamic systems, the corresponding dynamic control systems are also investigated. Some elementary properties are obtained.

The main purpose of this paper is to provide a framework for cross-dimensional linear systems. There are lots of problems remain for further investigation.
Particularly, the following problems are challenging and important:
\begin{itemize}
\item How to get general solution (trajectory) of cross-dimensional linear systems as $\A(t)$ is not dimension-bounded?
\item How to formulate the controllability and observability of continuous time cross-dimensional linear systems in general sense?
\item How to connect theoretical results with practical dimension varying systems, mentioned in the Introduction?
\item How to extend the results  in this paper to logical dynamic systems?
\end{itemize}

This paper can be considered as a follow up of \cite{chepr}. The next paper is concerning about the cross-dimensional nonlinear system. The overall purpose of this series papers is to build a theory on cross dimensional dynamic systems using dimension free matrix theory.

\end{document}